\documentclass[12pt,a4paper]{article}

\usepackage{pdflscape}
\usepackage{amsmath}
\usepackage{amssymb}
\usepackage{latexsym}
\usepackage{srcltx}
\usepackage{graphics}
\usepackage{color}
\usepackage{epsfig}
\usepackage{color}

\newcommand{\re}{{\mathbb R}}

\newcommand{\cA}{{\mathcal{A}}}

\newcommand{\cQ}{{\mathcal{Q}}}
\newcommand{\cP}{{\mathcal{P}}}

\newcommand{\cD}{{\mathcal{D}}}

\newcommand{\bx}{{\boldsymbol{x}}}
\newcommand{\by}{{\boldsymbol{y}}}
\newcommand{\bz}{{\boldsymbol{z}}}
\newcommand{\be}{{\boldsymbol{e}}}

\newcommand{\ba}{{\boldsymbol{a}}}
\newcommand{\bb}{{\boldsymbol{b}}}

\newcommand{\bv}{{\boldsymbol{v}}}
\newcommand{\bu}{{\boldsymbol{u}}}

\newcommand{\bs}{{\boldsymbol{s}}}

\newcommand{\vardot}{\mathord{\,\cdot\,}}

\newtheorem{theorem}{Theorem}
\newtheorem{prop}{Proposition}
\newtheorem{lemma}{Lemma}
\newtheorem{cor}{Corollary}
\newtheorem{remark}{Remark}

\date{}

\author{Vladimir Yu. Protasov 
\thanks{Moscow State University, Russia; University of L'Aquila, Italy, {e-mail: \tt\small
v-protassov@yandex.ru}} , 
Asiiat Musaeva
\thanks{Moscow State University, Russia; {e-mail: \tt\small
asya.musaeva2001@mail.ru
}} 
 }

\title{Second-order linear switching systems with arbitrary control sets: 
stability and invariant  norms
}

\begin{document}

\maketitle

\begin{abstract}

We show that the stability problem and the problem of constructing 
 Barabanov norms can be resolved for planar  linear switching systems 
 in an explicit form. This can be  done for  
every compact control set of $2 \times 2$ matrices.      
If the control set  does not contain a dominant matrix  with a real spectrum, 
then the invariant norm is always unique (up to a multiplier) and belongs to~$C^1$.  
Otherwise, there may be  infinitely many such norms, including non-smooth ones. 
All of them can be found and classified. In particular,  every 
symmetric convex body is a unit ball of the Barabanov norm of a suitable linear switching system. 

\bigskip

\noindent \textbf{Key words:} {\em linear switching systems, second-order, continuous-time, 
control set,  compact set of matrices, Lyapunov function, Barabanov norm, symmetric convex body, 
 stability, Frobenius norm 
}
\smallskip

\begin{flushright}
\noindent  \textbf{AMS 2020 subject classification} {\em 93D20, 52A10, 15A60}

\end{flushright}

\end{abstract}
\bigskip

\vspace{1cm}

\begin{center}

\large{\textbf{1. Introduction}}	
\end{center}
\bigskip

{\em Invariant} or {\em Barabanov} norm~$\re^d$ is defined for an arbitrary irreducible 
compact set~$\cA$ of~$d\times d$ matrices. It describes the asymptotic behaviour of trajectories 
 of a linear switching system, which is a linear differential equation with the matrix~$A(t)$ chosen from the set~$\cA$ at every 
 moment of time~$t$. The Barabanov norm not only defines the stability of the system but also gives a very precise information on the fastest  growth of its trajectories and, moreover, 
 finds those extremal trajectories.  Very  little is known about the properties of this norm: when it is unique, smooth, strictly convex? Apart from some trivial cases, 
there are  a  few examples  when the invariant norm can be  explicitly constructed. In this paper we 
prove that for $2\times 2$ matrices, this problem admits a complete solution. 
For every second-order linear switching systems, we classify all Barabanov norms, prove  a uniqueness  criterion (up to normalization),  and present  an algorithm for their construction. We will see that for a generic system,  this norm is unique and its sphere is a $C^1$ curve  described by one periodic trajectory. There is only one special case 
that provides a rich variety of invariant norms, including cases of non-smoothness and non-uniqueness. This is when 
  the set~$\cA$ contains a matrix with ``real dominance'' (all definitions are given below). 
 In this case, for  every norm in~$\re^2$, there is a linear switching system for which this norm is Barabanov. 
   In the geometrical interpretation this means that for every convex 
body~$G$ symmetric about the origin, there exists a compact set of matrices~$\cA$ whose 
Barabanov norm has the unit ball~$G$. 

It should be stressed, however, that the two-dimensional case is quite special and most of the results of this work are not extended, at least directly,  to higher dimensions. 
\medskip 

\textbf{Linear switching  systems}.  
We give all definitions in general dimension~$d$ and then focus on~$d=2$. 
 A {\em continuous-time linear switching system} corresponding to a given compact set~$\cA$
 of real $d\times d$ matrices is a linear ODE:  
 \begin{equation}\label{eq.main}
 \left\{
 \begin{array}{l}
  \bx(t) \ = \ A(t)\bx (t), \\
  \bx(0) \ = \ \bx_0, \\
  A(t) \in \cA, \ t\in [0, +\infty). 
 \end{array}
 \right. 
\end{equation}
The family~$\cA$ is called a {\em control set},  $\, A(\vardot) : \, [0, +\infty) \, \to \, \cA$ is an arbitrary measurable matrix function 
with values in~$\cA$ called {\em switching law}, the solution~$\bx(\vardot)\, : \, [0, +\infty)\, \to \, \re^d$
is a {\em trajectory} of the system. We will identify system~(\ref{eq.main}) with its control set~$\cA$ since it is uniquely defined by this set.  The maximal asymptotic growth of trajectories, 
over all switching laws~$A(\vardot)$ and all initial points~$\bx_0$, 
is characterized by the Lyapunov exponent $\sigma(\cA)$. This is the 
infimum of numbers~$\sigma$ such that for every trajectory~$\bx(\vardot)$, 
there is a constant~$C = C(\bx)$ for which \linebreak $\|\bx(t)\| \ \le \ C\, e^{\sigma t}, \ 
t \ge 0$.
\medskip

\textbf{Stability and the Barabanov norms}.  
The stability of linear switching systems has been systematically studied  
since 70th, see, for instance~\cite{B, BT79, BT80, MNV84, MP86, O}. 
The system is (asymptotically) stable if all its trajectories tend to zero. 
We only deal with this type of stability and usually omit the word ``asymptotically''.  
Obviously, if~$\sigma(\cA) < 0$, then the system is stable. The converse is less obvious, but also true~\cite{MP89}. The Lyapunov exponent does not change after replacing the control set~$\cA$ by its convex hull~${\rm co}(\cA)$. 
A {\em shift} of the set~$\cA$ by a number~$\alpha$ is $\cA - \alpha I \, = \, 
\{A-\alpha I: \ A \in \cA\}$, where~$I$ is the identity matrix. 
Every trajectory of the system~$\cA$ can be written as~$\bx(t) = e^{\,\alpha t} \bx_{\alpha}(t)$, where~$\bx_{\alpha}$ is the trajectory of the shifted family~$\cA - \alpha I$ corresponding 
to the same switching law and to the same~$\bx_0$. Therefore, $\sigma(\cA - \alpha I) \, = \, 
\sigma(\cA) - \alpha$. In particular, after the shift by~$\alpha = \sigma(\cA)$ 
the Lyapunov exponent becomes zero. That is why we often assume that 
 $\sigma(\cA) = 0$, i.e., that we passed from the original family~$\cA$ to~$\cA - \sigma I$. 
 We never consider the matrix~$\sigma I$ if it belongs to~$\cA$. 
 Indeed, after the shift by~$\sigma$ it becomes the zero matrix and its every trajectory is constant. 
 Hence, it does not influence the other trajectories and does not change the Lyapunov exponent. In particular, when we say that in case~$\sigma(\cA) = 0$ the system 
 does not contain degenerate matrices, we do not count the zero matrix. 
 \smallskip 
 
 A norm~$f(\bx)$ in~$\re^d$ is called {\em invariant norm} or {\em Barabanov norm}
 for system~(\ref{eq.main}), if every trajectory~$\bx(\vardot)$
 satisfies~$f(\bx(t))\, \le \, e^{\sigma t} f(\bx_0), \, t\ge 0$, 
  and for every~$\bx_0 \in \re^d$, there exists an~{\em extremal trajectory}~$\bar \bx(\vardot)$ of the system with  the control 
 set~${\rm co}(\cA)$
 such that~$\bar \bx(0) = \bx_0$ and 
 $f(\bar \bx(t))\, = \, e^{\sigma t} f(\bx_0), \, t\ge 0$. Thus, the extremal trajectories pass through every point and they are trajectories of the maximal growth. 
 If~$\sigma(\cA) = 0$, which can be assumed after the shift~$\cA \, \mapsto \, \cA - \sigma I$, 
 then the definition of the Barabanov norm can be formulated as follows: 
 for every trajectory~$\bx(\vardot)$, the norm~$f(\bx(t))$ is non-increasing in~$t$
 and for every point~$\bx_0$, there exists an extremal trajectory~$\bar \bx(\vardot)$ 
 (corresponding to the control 
 set~${\rm co}(\cA)$)
 starting at~$\bx_0$ and such that~$f(\bar \bx(t)) \equiv f(\bx_0)$. 
 
In case~$\sigma(\cA)= 0$ there exists an equivalent 
geometrical definition of invariant norm. Let~$G$ be the unit ball of a norm~$f$
and~$S = \partial G$ be its unit sphere. 
 Then~$f$ is Barabanov  if and only if for every~$\bx \in S$, we have: 
 1) for each~$A \in \cA$, the vector~$A\bx$ starting at~$\bx$ is either 
 tangent to~$S$ or directed inside~$S$ and 2) the set of  tangent vectors~$A\bx$
 is nonempty. Here ``directed inside'' means that~$\bx + \lambda A\bx \in {\rm int}\, G$
 for small~$\lambda > 0$; ``is tangent to~$S$'' means that the distance 
 from the point~$\bx + \lambda A\bx$ to~$S$ is~$o(\lambda)$ as~$\, \lambda \to +0$. 
 In particular, the zero vector is tangent to~$S$. 
 
 For all control sets~$\cA - \alpha I, \, \alpha \in \re$, the Barabanov norm is the same.  Therefore, 
 it suffices to take~$\alpha = \sigma$, i.e., to consider those norms only for the case~$\sigma(\cA) = 0$. 
The sets~$\cA$ and~${\rm co}(\cA)$ have the same Barabanov norm. 

A linear switching system  is called {\em irreducible} if the matrices of its control 
set~$\cA$ do not share a common invariant 
subspace.   It was proved in~\cite{B} that a compact irreducible system
always possesses an invariant norm. This norm may not be unique. Here and below we mean the uniqueness up to multiplication  by a constant. 

\smallskip 

\medskip

\textbf{Continuous-time vs discrete-time.} 
Barabanov norms have been intensively studied in the literature, not only 
for system~(\ref{eq.main}) but also for discrete-time systems, when 
the ODE is replaced by the difference equation~$\bx(k+1) = A(k)\bx(k)$. 
For such systems, the invariant norms are well understood, see, for instance~\cite{C, GZ1, Mo}. 
For most of  discrete-time systems 
with a finite set of matrices, Barabanov's norms are unique~\cite{P} and have a simple structure: they are either piecewise-linear or 
piecewise-quadratic~\cite{GZ2, MeP}. The situation  for continuous-time systems~(\ref{eq.main})
is  more complicated, see~\cite{CGM, GCM, GO, TM} and references therein. As a rule, 
nothing is known on the invariant norms except for their existence. 
In the works~\cite{CGM, GCM} several important 
properties of Barabanov's norms in dimensions~$d=2, 3$ have been obtained under
additional assumptions.   In some works the invariant 
norms have been approximated by simpler norms, such as piecewise-linear, etc.~\cite{BM96, BM99, GLP, OIGH93}.

 We are going to see that 
in the two-dimensional case the invariant norms can be found, constructed  and  completely classified for every compact set of matrices. 
\medskip

\textbf{The case of $2\times 2$ matrices.}  Stability and invariant norms in dimension~$d=2$
have been studied in great detail~\cite{BB, BBMZ, GSM, GTI, HM, TM, S}. While  for discrete-time systems, the 2D case is not much simpler than the higher dimensions~\cite{C, GZ1},   
in continuous-time the planar case is indeed special. 
For example,  the stability of finite sets of~$2\times 2$ matrices can efficiently 
be decided. For sets of two matrices, this was done in~\cite{BBM,  TM, TA} by the search 
of the ``worst trajectory'' and by the tools of Lie algebra~\cite{Ma}, see also~\cite{BB, HM, S, Y1, Y2}.  
For arbitrary finite sets of matrices, this was done in~\cite{GTI, Mu1}. 

In all the aforementioned works 
the finiteness of the control set~$\cA$ is significant. 
All suggested methods used an exhaustion of~$\cA$.  The question arises: can the stability problem
 be solved explicitly for an arbitrary compact set~$\cA$? 
In Theorem~\ref{th.30} 
we prove a stability criterion  for arbitrary compact family of~$2\times 2$ matrices. The invariant norm issue for planar continuous-time systems has also been  addressed only 
for finite matrix sets~\cite{CGM, GCM, TM, Mu1}. We discuss this issue in~Section~6.  
In Section~4 of  this paper we obtain a complete solution. Moreover, we show that the set 
of all possible invariant norms is wide: every norm in~$\re^2$ is Barabanov for a suitable set~$\cA$.  
\medskip

\textbf{The main results.} We solve the stability problem and find the Barabanov norms 
for arbitrary compact set of $2\times 2$ matrices. The stability is addressed in 
Section~3, here we describe the results on the invariant norms. 

 A matrix~$A \in {\rm co}(\cA)$ (different from~$\sigma I$), is called 
{\em dominant} if the largest real part of its eigenvalues is equal to~$\sigma(\cA)$.
In this case the trajectory~$\bx(t) \, = \, e^{tA}\bx_0$ generated by the 
stationary switching law~$A(t) \equiv A$, has the maximal growth 
unless~$\bx_0$ belongs to an invariant subspace of~$A$ corresponding to 
eigenvalues with smaller real parts. 
For $2\times 2$ matrices, there are two possible cases of dominance:   
 {\em real}, if $A$ has real eigenvalues, and  {\em complex} otherwise. 
 The matrix~$\sigma I$,  if it belongs to~$\cA$, is not considered 
 and is not counted as dominant. 
\smallskip

Our main result  asserts that for an arbitrary compact set 
of $2\times 2$ matrices, there are three possible cases for the Barabanov norm:  
  \smallskip

1) If the system does not have a matrix with real dominance, then 
its Barabanov norm is unique and its  unit sphere~$S$
coincides with one period of 
a certain periodic trajectory of the system~$\cA - \sigma I$. 
This trajectory can be explicitly computed. 
\smallskip  

2) If the system has a matrix~$A$ with complex dominance, then 
its Barabanov norm is quadratic, the unit sphere is an ellipse
described by the trajectory of the matrix~$A - \sigma I$. 
\smallskip  

3) If the system has real dominance, then there may be infinitely many 
Barabanov norms. Their complete variety    is given in Theorem~\ref{th.40}.
The unit sphere of the Barabanov norm is never described by one dominant trajectory. 
Every convex curve symmetric about the origin can be a unit sphere of the Barabanov norm 
of a suitable system with  real dominance.   
\medskip 

Note that the dominance property (real or complex) can be identified  without involving 
the convex hull of~$\cA$. Theorem~\ref{th.35} proven in Section~4 asserts that 
if ${\rm co}(\cA)$ contains a dominant matrix, then so does~$\cA$. 

\medskip 

\textbf{Other results.} The stability criterion given by Theorem~\ref{th.30} 
makes it possible to decide between the cases~$\sigma(\cA) < 0$ and~$\sigma(\cA) \ge 0$. 
Applying bisection in the shift of the set~$\cA$, one can compute the Lyapunov exponent
~$\sigma(\cA)$ with an arbitrary precision. After this the invariant norm can be constructed algorithmically. The construction is rather simple in the case when the set~$\cA$ has no 
real dominance (Section~2). Otherwise,  it is more technically difficult due to 
a wide set of invariant norms. The algorithm is given in subsection~4.3.   
As a simple corollary we obtain that every 2D system has an extremal trajectory
generated by a periodic switching law.

The uniqueness issue for the invariant norms is solved in Section~4. 
Theorem~\ref{th.50} gives the corresponding criterion. 

Section~7 deals with applications.  
The first one concerns the stability of a matrix ball~$\cA \, = \, \{A: \ \|A-A_0\|_F \le r\}$, where~$\|\cdot \|_F$ is the Frobenius norm. The second problem is the stability 
of the system~(\ref{eq.main}) with the control set given approximately, with possible noise.

\medskip 

\textbf{Novelty.} Most of the literature on the stability of second-order linear switching systems 
 deal with the case of two matrices. Finite matrix sets~$\cA = \{A_1, \ldots , A_m\}$ 
 for arbitrary~$m\ge 2$ are considered in~\cite{GTI, Mu1}. The methods derived in those works demonstrate that 
 their complexity growth dramatically with the number of matrices. 
 For example, the algorithm for deciding stability presented in~\cite{Mu1} requires  about~$m^2/2$ cases. Even an approximate solution  using the  popular CQLF method 
 requires solving a system of~$m$ linear matrix inequalities, which becomes hard for large~$m$. 
 Therefore, replacing a generic convex set of matrices~$\cA$ by a polytopic set 
 is not efficient even for an approximate solution. We suggest a new method 
 which works equally well for generic sets.  As a special case we consider matrix balls and apply it to the 
 stability problem with noisy data. Those problems could not  be solved by the known methods. 
 
  The construction of the Barabanov norm is realized in an explicit form. We also provide a complete classification  of those norms and solve the uniqueness problem. We show that in the general case, when the set~$\cA$ has no real dominance, the Barabanov norm is unique 
  (up to multiplication) and is defined by a smooth periodic trajectory. The case of real dominance, being quite special, generates a rich variety of Barabanov norms and 
  admits non-uniqueness. Apart from the
  theoretical interest those results lead to sharp estimates of the fastest growth of trajectories.

\medskip 

\textbf{Notation.} If the converse is not stated, we always deal with 
 an irreducible compact set~$\cA$ of $2\times 2$ matrices. The unit ball of the 
 Barabanov norm is denoted by~$G$ and the unit sphere (actually, a closed convex curve 
 symmetric about the origin) by~$S$. We use the standard terminology for 
 lines of support for~$S$ and for left and right tangent lines. 
 Denote by~$\alpha_{\bx}$ the {\em left tangent ray} to~$S$ at the point~$\bx\in S$. 
 This is the limit of rays $\ell_{\by} = \{\bx + \lambda (\by - \bx): \ \lambda \ge 0\}$
 as~$\by$ tends to~$\bx$ along~$S$ from the left, i.e.,  
 the arc $\bx \by$ has positive direction. Similarly we define the right tangent 
 ray~$\beta_{\bx}$. If those rays are not collinear, then~$\bx$ is called {\rm a corner point}. 
 The set of corner points is finite or countable. 
 For all other points of~$S$,  the line of support is unique and 
 coincides with the tangent (and hence, the left and right tangents coincide).

 We call a matrix Hurwitz if all its  eigenvalues have negative real parts. 
 
 As usual, the positive direction of a rotation ({\em positive rotation}) is counterclockwise, the negative one is clockwise. The oriented angle between vectors 
 $\ba$ and $\bb$ is the angle of positive rotation that maps the direction of~$\ba$ to 
 the direction of~$\bb$. The same with the angle between rays. 
\begin{remark}\label{r.5}
{\em Most of the literature on linear switching systems deals with finite sets of 
matrices~$\cA = \{A_1, \ldots , A_m\}$. In this case the function~$A(t)$
has finitely many values with ``switches'' between them. For arbitrary compact sets~$\cA$, 
the values of~$A(t)$ may change continuously, in which case it does not  have real switches. 
Nevertheless, in many works the term ``switching'' is used  for arbitrary control sets as well. 
We follow this terminology. We also refer to the matrix set~$\cA$ as 
``control set'', although in the optimal control this term usually has a different meaning. We do this for the sake of simplicity. 
}
\end{remark}

\vspace{1cm}

\begin{center}
\large{\textbf{2. The Barabanov norms for  systems without real dominance}} 
\end{center}
\bigskip 

In this section we consider sets~$\cA$ for which none of the matrices from~${\rm co}(\cA)$ 
has real dominance, i.e., 
has an eigenvalue equal to the Lyapunov exponent~$\sigma(\cA)$ of the system~(\ref{eq.main}). 
 This means 
that none of the stationary switching laws~$A(t) \equiv A$
generated by one matrix~$A\in {\rm co}(\cA)$ with  a real leading eigenvalue 
provides  the maximal asymptotic growth. 
We are going to see that in this case there exists a unique Barabanov norm. 
This norm is 
generated by a periodic trajectory. This unique extremal trajectory may have infinitely many switching points and can be  obtained explicitly as a solution of a certain linear differential equation. On the other hand, there is one particular case when this extremal trajectory has no 
switching points at all. This is the case of complex dominance, with which 
we begin our analysis. In this case the unique Barabanov norm is quadratic and its unit ball is an ellipse. 

\medskip

\newpage

\begin{center}
\textbf{2.1. Systems with complex dominance} 
\end{center}
\bigskip 

As usual, a quadratic norm is~$f(\bx) \, = \, \sqrt{\bx^TM\bx}$, where~$M$ 
 is a positive definite matrix.
We begin with the following simple  observation: 
\begin{prop}\label{p.10}
If $\sigma(\cA) = 0$ and the convex hull~${\rm co}(\cA)$
contains a matrix with a  purely imaginary spectrum, then 
all such matrices from~${\rm co}(\cA)$ are proportional, at least one of them belongs to~$\cA$
  and a unique Barabanov norm of~$\cA$
is quadratic.    
\end{prop}
{\tt Proof}. Every $2\times 2$-matrix with an imaginary spectrum
is similar to a rotation by~$\frac{\pi}{2}$ multiplied by a scalar. 
Passing to the corresponding basis we obtain a matrix~$A$, for which 
all trajectories~$\dot \bx = A\bx$ are concentric circles. A Barabanov 
norm of~$\cA$ is non-increasing along those  trajectories, 
hence, it is constant since the trajectories are periodic. 
Thus, the Barabanov norm is proportional to the Euclidean one and, therefore, it is 
unique up to normalization.  
If there is another operator~$A' \in {\rm co}(\cA)$ 
with an imaginary spectrum, then it also generates a Barabanov norm, 
which is Euclidean in another basis. The uniqueness implies that 
these two bases are similar by means of an orthogonal transformation. 
Hence, $A'$ is also proportional to a rotation by~$\frac{\pi}{2}$. 
Thus, all matrices from~${\rm co}(\cA)$ with imaginary spectra are proportional. 
Assume that none of them belongs to~$\cA$. In this case, the matrix~$A$
is a convex combination of some~$A_1, \ldots ,A_N \in \cA$. 
For every~$\bx$ from the unit circle~$S$, all the vectors~$A_1\bx, \ldots , A_N\bx$ 
starting at~$\bx $ are either 
tangent to $S$ or directed inside~$S$. Their convex combination~$A\bx$ is tangent to~$S$, hence, 
all $A_i\bx$ are also tangent. Thus, $A_i\bx \perp \bx$  for all~$\bx \in \re^2$, consequently, $A_i$ is proportional to the rotation by~$\frac{\pi}{2}$, which is a contradiction.  

  {\hfill $\Box$}
\medskip 

Thus, in the case of complex dominance, the unique Barabanov norm is affinely similar to the Euclidean one, 
its invariant sphere is an ellipse~$S=\{\bx \in \re^2: \, \bx^TM\bx = 1\}$, where~$M \succ 0$.  It is easy to characterize the families 
with complex dominance and to find~$S$. 

\begin{theorem}\label{th.15}
A family~$\cA$ has  complex dominance if and only if there exists a matrix~$A_0 \in \cA$
with eigenvalues~$\alpha \pm \beta i$ such that~$\beta \ne 0$ and the system 
of linear matrix inequalities~$A^TM + MA \, \preceq \, \alpha I, \, A \in \cA$,  possesses a solution~$M\succ 0$. In this case~$\sigma(\cA) \, = \, \alpha\, , \, A_0$ is dominant, 
the Barabanov norm is unique  and is equal to~$f(\bx) = \sqrt{\bx^TM\bx}$.  
\end{theorem}
{\tt Proof}. After the shift of the system~$\cA \mapsto \cA - \alpha I$
it can be assumed that~$\alpha = 0$. It is well-known that if a matrix~$M\succ 0$
satisfies the inequalities~$A^TM + MA \, \preceq \, 0, \, A \in \cA$, then  the norm~$f(\bx) = \sqrt{\bx^TM\bx}$ is non-increasing
 along any trajectory~\cite{L}. Therefore,~$\sigma(\cA) \le 0$.   On the other hand,~$\sigma(A_0) = 0$, consequently~$\sigma(\cA) \ge 0$. Thus, $\sigma(\cA) \ge 0$ and therefore, 
 the system has complex dominance. Conversely, if the system has complex dominance, then 
 Proposition~\ref{p.10} implies the existence of a quadratic 
 Barabanov norm~$f$. This norm is non-increasing along every trajectory, 
 hence, its matrix~$M$ has to satisfy the system~$A^TM + MA \, \preceq \, \alpha I, 
 \, A \in \cA$.

  {\hfill $\Box$}
\medskip

\medskip

\newpage

\begin{center}
\textbf{2.2. Systems without real  dominance} 
\end{center}
\bigskip

A point~$\bx\in \re^2$ is called {\em feasible} if   
there is no matrix~$A \in \cA$ and number~$\lambda \ge 0$ such that~$A\bx = \lambda \bx$. 
Thus, $\bx$ is not an eigenvector of a matrix from~$\cA$ with a nonnegative eigenvalue. 
Consider the set of images~$\{A\bx\, : \ A \in \cA\}$. 
If it contains a vector~$\by$ such that 
the  (oriented) angle between~$\bx$ and~$\by$ belongs to the 
interval~$(0, \pi)$, then we take the vector for which this angle is 
minimal (it  exists under our assumptions) and, among them, 
the vector~$\by$ with the maximal length.
The matrix~$A \in \cA$ for which~$A\bx = \by$ 
will be denoted by~$A^{(l)}_{\bx}$ and called the {\em leading left matrix}. 
If there are several such matrices, then we take any of them.  
This matrix is well-defined whenever 
$\bx$ is feasible  and there is~$A \in \cA$ 
such that the angle between~$\bx$ and the image~$A\bx$ 
belongs to~$(0, \pi)$. Similarly we define 
the {\em leading right matrix~$A^{(r)}_{\bx}$}, for which the image 
makes the maximal angle with~$\bx$ on the interval~$(-\pi, 0)$. 
We usually omit the subscript~$\bx$ keeping in mind that~$A^{(l)}$ and~$A^{(r)}$
depend on~$\bx$. The rays starting at the point~$\bx$
and directed along~$A^{(l)}\bx$ and~$A^{(r)}\bx$ 
will be denoted by~$a_{\bx}$ and~$\, b_{\bx},\, $  respectively.
These are the {\em leading left} and the {\em leading right} directions. 
They correspond, respectively, to the positive and negative directions of the rotation about the 
origin. 

For infeasible~$\bx$, the leading directions are not defined, 
for feasible~$\bx$ at least one of them is defined.
Otherwise,~$A\bx$ is collinear to~$\bx$
 for all~$A\in \cA$, hence~$\cA$ is reducible. 
\begin{lemma}\label{l.5}
If the left leading direction~$a_{\bx}$ is defined at a feasible point~$\bx$, then 
it is defined in a neighbourhood of~$\bx$ and continuously depends on the point.
\end{lemma}
{\tt Proof}. A point~$\bx$ is feasible if and only  if  the compact set~$\cA_{\bx}\, = \, \{A\bx : \ A\in {\rm co} (\cA)\}$
does not intersect the ray~$\{\lambda \bx: \ \lambda \ge 0\}$. Since~$\cA_{\bx}$ 
continuously (in the Hausdorff distance) depends on~$\bx$, then all close points are also 
feasible and  the directions of vectors from~$\cA_{\bx}$  forming maximal and minimal angles with~$\bx$ continuously depends on~$\bx$.

  {\hfill $\Box$}
\medskip

The {\em leading left trajectory} is the solution of ODE~$\dot \bx (t) = A^{(l)}(t)\bx(t)$
on a segment~$[t_0, t_1]$ with some initial condition~$\bx(t_0) = \bx_0 \ne 0$, where 
~$A^{(l)}(t)$ is the operator~$A^{(l)}$ associated to~$\bx(t)$. 
It is assumed that all~$A^{(l)}(t)$ are well-defined for all~$t\in [t_0, t_1]$.  
This trajectory goes around the origin in the positive direction 
(counterclockwise). 
The same is for the leading right trajectory~$\dot \bx (t) = A^{(r)}(t)\bx(t)$. 
Lemma~\ref{l.5} implies that both the leading trajectories are~$C^1$. 

\begin{theorem}\label{th.20}
Suppose~$\sigma(\cA) =0$ and the system does not have  real dominance 
(i.e., ${\rm co}(\cA)$ does not contain 
degenerate matrices); then the Barabanov norm is unique and its unit sphere~$S$ 
is described by one period of a periodic leading trajectory, either left or right. 

Moreover, if the system does not have  complex dominance, 
then it possesses a unique, up to normalization,  leading trajectory. 
 If this  is a left trajectory, then  at every point~$\bx \in S$, the 
 ray~$a_{\bx}$ is tangent to~$S$ and~$b_{\bx}$ is either  directed inside~$S$ or 
does not exist.  In  case of the right trajectory, it is vice versa. 
\end{theorem}
\begin{remark}\label{r.35}
{\em Theorem~\ref{th.20} can be formulated without involving the 
convex hull of~$\cA$.  Theorem~\ref{th.35} in Section~4 asserts that 
if ${\rm co}(\cA)$ contains a dominant matrix, then so does~$\cA$, unless~$\cA$
has complex dominance. }
\end{remark}

{\tt Proof of Theorem~\ref{th.20}} is realized in four steps:  
\smallskip 

1. {\em For every~$\bx \in S$, 
either one of the vectors~$a_{\bx}, b_{\bx}$ is not defined or 
the oriented angle between~$a_{\bx}$ and~$\, b_{\bx}$
is less than~$\pi$}. 

Otherwise, the ray~$\{\lambda \bx:  \ \lambda \ge 1\}$
 intersects 
the segment connecting points~$\bx + A^{(l)}\bx$ and~$\bx + A^{(r)}\bx$. This means that 
there exists a convex combination~$s A^{(l)}\bx \, + \, (1-s)A^{(r)}\bx, \, s\in [0,1]$,  
which is  directed along~$\bx$. Then the matrix~$A \, = \, 
s A^{(l)}\, + \, (1-s)A^{(r)}$ has a real nonnegative 
eigenvalue. If it is zero, then~$A$ is 
degenerate; if it is positive, then~$\sigma(\cA)  > 0$. The contradiction completes the proof
of step~1. 
\smallskip 

2. {\em The  set of points~$\bx \in S$ at which both~$a_{\bx}$ and~$b_{\bx}$
are tangent to~$S$ is either empty or finite. We call such points~{\em switching}. 
On every open arc~$\gamma \subset S$ 
connecting neighbouring switching points, either all~$a_{\bx}, \, \bx \in \gamma$, are 
tangents and all~$b_{\bx}$ are not, or vice versa. This arc is a~$C^1$ curve.} 

According to Step~1, the angle between~$a_{\bx}$ and~$b_{\bx}$
is less than~$\pi$. Assume there is a sequence of points~$\bx_i \in S$
such that the angle between~$a_{\bx_i}$ and~$b_{\bx_i}$ tends to~$\pi$. By the compactness, passing to a subsequence 
it can be assumed that~$\bx_i$ tends to some~$\bar \bx\in S$ and 
the matrices~$A^{(l)}_{\bx_i}, A^{(r)}_{\bx_i}$ converge  to some matrices~$\bar A_l, \bar A_r \in {\rm co}(\cA)$, respectively.  
Therefore~$\lim_{i \to \infty} A^{(l)}\bx_i \, = \, \bar A_l\bar \bx$
and~$\lim_{i \to \infty} A^{(r)}\bx_i \, = \, \bar A_r\bar \bx$.    
Note that the vectors~$\bar A_l\bar \bx$ and~$\bar A_r\bar \bx$ 
are nonzero (since~${\rm co}(\cA)$ does not contain degenerate matrices) and   
the angle between them is equal to~$\pi$. 
Hence, at least one of the rays~$a_{\bar \bx}, b_{\bar \bx}$  is not defined. 
Otherwise, they are directed along~$\bar A_l\bar \bx, \bar A_r\bar \bx$, respectively, and form the angle~$\pi$. 
Assume~$a_{\bar \bx}$ is not defined. 
In  this case~$\bar A_l \bar x = \lambda \bar \bx$ for some~$\lambda \in \re$. 
 If~$\lambda = 0$, then~$\bar A_l$ is degenerate; 
if~$\lambda > 0$, then~$\sigma(\bar A_l) > 0$. Both those cases are impossible, hence~$\lambda < 0$. Thus, the vector $\bar A_l \bar \bx$ is directed along~$- \bar \bx$. 
Since it makes the angle~$\pi$ with the vector~$\bar A_r \bar \bx$, we see 
that~$\bar A_r \bar \bx$
is directed along~$\bar \bx$, i.e.,~$\bar A_r\bar \bx = \mu \bx$
for some~$\mu > 0$. Hence, $\sigma(\bar A_r) \ge \mu > 0$, which is a contradiction.  Thus, there is no sequence~$\bx_i$ along which 
the angle between~$a_{\bx_i}$ and~$b_{\bx_i}$ tends to~$\pi$. Therefore, 
there is~$\varepsilon > 0$
 such that 
the  angle between $a_{\bx}$ and~$b_{\bx}$ does not exceed~$\pi - \varepsilon$ at all points~$\bx \in S$.

If at some point~$\bx$ both $a_{\bx}$ and~$b_{\bx}$  touch~$S$, 
then~$\bx$ is a corner point of~$S$ such that the adjacent angle between one-side tangents 
is at least~$\varepsilon$. Since the sum of all  adjacent angles at the corner points of
a convex curve does not exceed~$2\pi$, it follows that the total number of corner points~$\bx$
does not exceed~$2\pi / \varepsilon$.  On every arc between neighbouring points, 
the set of points~$\bx$ for which~$a_{\bx}$ touches~$S$ is closed.  The same is true for~$b_{\bx}$. 
If a connected set is split to two closed sets, then one of them is empty. 
Thus, every arc has one direction  of tangent rays ($a_{\bx}$ or~$b_{\bx}$). 
Since by Lemma~\ref{l.5} those rays depend continuously on~$\bx$, it follows that the arc is~$C^1$. This lemma is applicable since under our assumptions, all points~$\bx\ne 0$
are feasible. 
 \smallskip 

Thus,~$S$ admits a finite partition by switching points, 
each interval has one direction of tangents: either~$a_{\bx}$ or~$b_{\bx}$. 
 Now we are going to prove that 
this partition is actually trivial and all tangents  have the same direction on the entire 
curve~$S$. 
 \smallskip 

3. {\em Either all the rays~$a_{\bx}, \, \bx \in S$, touch~$S$ and all~$b_{\bx}$
do not, or vice versa. } 


\smallskip 

Assume the contrary:  there are several  arcs of the partition
by switching points. 
After concatenation of arcs corresponding to the same direction
($a_{\bx}$ or~$b_{\bx}$) it can be assumed 
that neighbouring arcs have opposite directions of tangent vectors. 
Going along~$S$ in positive direction, we 
meet two neighbouring arcs such that the first 
arc~$\bx_1 \bx_2$  (the enumeration is counterclockwise)
has tangents~$a_{\bx}$ and the next arc~$\bx_2\bx_3$ 
has tangents~$b_{\bx}$. 

Since each arc is~$C^1$, it follows that the right tangents (that are opposite to~$a_{\bx}$)
tend to the right tangent at~$\bx_2$ as~$\bx\to \bx_2$ along the arc~$\bx_1\bx_2$. On the other hand,  
$a_{\bx}$ continuously depends on~$\bx$ (Lemma~\ref{l.5}), therefore, 
$a_{\bx}$ tends to~$a_{\bx_2}$ as~$\bx\to \bx_2$ along the arc~$\bx_1\bx_2$. 
Thus, $- a_{\bx_2}$ is the right tangent to~$S$ at~$\bx_2$. 
Similarly, $-b_{\bx_2}$ is the left tangent at~$\bx_2$. 
By the convexity, the oriented angle between the left and the right tangents 
belongs to~$(0, \pi]$. Hence, so does the angle between~$b_{\bx_2}$ and~$a_{\bx_2}$. 
We see that the oriented angle between~$a_{\bx_2}$ and~$b_{\bx_2}$ is bigger than or equal to~$\pi$
which contradicts  Step 1.  

 \smallskip

4. {\em One of the two trajectories~$\dot \bx = A^{(l)}\bx$ or~$\dot \bx = A^{(r)}\bx$ starting 
from an arbitrary point of~$S$
is periodic and coincides with~$S$. 
The other one is either non defined or converges  to zero}.

This follows directly from Step~3: if, say, all~$a_{\bx}$ are tangent to~$S$, then $S$ 
is a periodic trajectory of the equation~$\dot \bx = A^{(l)}(t)\bx$. 

  {\hfill $\Box$}
\medskip 

Now we can make the main conclusion for 
general systems, without the condition~${\sigma (\cA) = 0}$. 

\begin{cor}\label{c.10}
A system without real dominance always possesses a unique Barabanov norm. The unit 
sphere~$S$ of this
norm is a~$C^1$ curve. 

If, in additional, it does not have  complex dominance, then 
it possesses a unique extremal  trajectory~$\bar \bx(\vardot
)$. It corresponds to a periodic switching law and 
the curve~$e^{-\sigma t}\bar \bx(t), \, t\in [0, 2T)$, where~$T$ is the period, 
describes~$S$. 

\end{cor}
{\tt Proof}. By the corresponding shift~$\cA \mapsto \cA - \sigma I$
we pass to the case~$\sigma(\cA) = 0$ and imply Theorem~\ref{c.20}. 
The uniqueness follows immediately, the~$C^1$ regularity follows from 
the continuity of~$a_{\bx}$ and~$b_{\bx}$ with respect to~$\bx$. 

  {\hfill $\Box$}
\medskip 

Theorem~\ref{th.20} immediately implies that the fastest asymptotic growth 
of trajectories of an irreducible 2D systems is always achieved on periodic switching laws.  
\begin{cor}\label{c.12}
A second-order linear switching system with a convex control set~$\cA$ and with~${\sigma(\cA) = 0}$ always possesses a  periodic extremal trajectory. 
\end{cor}
{\tt Proof}. If a system does not have  real dominance, then 
the statement follows from Theorem~\ref{th.20}. 
Otherwise, there exists a matrix~$A_0\in \cA$ with real dominance, i.e., 
$A_0$ is degenerate. Taking an arbitrary point~$\bx_0 \ne 0$ from the kernel of~$A_0$, 
we obtain a periodic trajectory~$\bx(t)\, \equiv \, \bx_0$.  

  {\hfill $\Box$}
\medskip 
 
 Thus, the fastest asymptotic growth of trajectories of a planar system
 is always achieved at a  periodic switching law. For control sets of two matrices,  
 Corollary~\ref{c.12} was well-known, see~\cite{BB, BBM, GTI, TA, Y2}, where the periodic extremal trajectories 
 were constructed in an explicit form.

For discrete-time linear switching systems, this is 
 not true. It was a long-standing Lagarias-Wang finiteness  conjecture (1995)
eventually disproved in~\cite{BTV} for planar systems with two-matrix switching sets. On the other hand, 
it is likely that all ``generic'' discrete-time system (of an arbitrary order~$d$)  do possess 
this property, see~\cite{GP, GZ2}. As for the continuous-time systems, this is still unknown whether 
Corollary~\ref{c.12} holds for arbitrary dimensions~$d$~? 
Or, at least, is it true for all ``generic'' in some sense linear switching systems? 

\begin{remark}\label{r.45}
{\em Theorem~\ref{th.20} can be put into practice by finding the 
left and right leading trajectories of~$\cA$. Finding~$a_{\bx}$ and~$b_{\bx}$
for each~$\bx$ we write the differential equations~$\dot \bx = A^{(l)}\bx$
and~$\dot \bx = A^{(r)}\bx$ with an arbitrary initial point~$\bx_0$, from which we 
compute the leading trajectories. If~$\sigma(\cA) = 0$ and there is no real dominance, then 
one of those trajectories is periodic and defines the unit sphere~$S$ of the Barabanov norm}. 
\end{remark}

\begin{center}
\textbf{2.3. Positive systems always have real dominance} 
\end{center}
\bigskip

A system is called positive if~$\cA$ consists of {\em Metzler matrices}, 
whose off-diagonal elements are all nonnegative. All trajectories of a positive system 
started in the positive quadrant~$\re^2_+$ never leave this quadrant~\cite{GSM, TM}. 
  A positive system cannot have periodic trajectories making rotations around the origin.  
Otherwise, such a  trajectory  contains some positive vector~$\bx (t_1)$ and hence 
$\bx (t)$ is also nonnegative for all~$t\, > \, t_1$. Hence, this trajectory will never 
come to the point~$-\bx(t_1)$. 
Thus, we have 
\begin{cor}\label{c.20}
A positive planar system always has 
real dominance.  
\end{cor}
{\tt Proof.} After a proper shift we assume that~$\sigma(\cA)=0$. 
Theorem~\ref{th.20} implies that either~${\rm co}(\cA)$ contains a 
degenerate  matrix, which is dominant, or it has a periodic trajectory. 
The latter is impossible for a positive system.

  {\hfill $\Box$}
\medskip 

Corollary~\ref{c.20} immediately implies Theorem 3.3 from~\cite{GSM}, which states that
 {\em a positive planar system with all  stationary trajectories converging to zero is stable}. Indeed, this system has a dominant trajectory which also has to converge to zero. This means that the system is stable.

\medskip

\begin{center}
\large{\textbf{3. Stability}}
\end{center}


\bigskip 

We call a leading trajectory~$\bx (\vardot)$ {\em round} if 
the vector~$\bx(t)$ makes a rotation about the origin. This means that
there exists~$T>0$ such that 
$\bx(T)  \, = \,  -\, \lambda \, \bx(0)$ for some~$\lambda > 0, \, T > 0$, and 
the switching law~$A(\vardot)$ is periodic with the minimal period~$T$. 
If~$\lambda = 1$, then the trajectory is periodic;  
if~$\lambda < 1$ ($\lambda > 1$), then the trajectory is called~{\em decreasing} (respectively, {\em increasing}). Every system can have 
at most two round leading trajectories (left and right) but may not have them at all.

\begin{theorem}\label{th.30}
 An irreducible  family~$\cA$ is  stable if and only if ${\rm co} (\cA)$ does not contain 
matrices with a nonnegative  eigenvalue and the both round leading trajectories 
are  decreasing, whenever they exist. In this case 
an arbitrary point~$\bx_0 \ne 0$ can be chosen as a starting point of these trajectories. 
\end{theorem}
{\tt Proof.}  {\em (Necessity)}.  If 
a leading trajectory is non-decreasing, then 
it does not converge to zero, so, the system is unstable. 
If some of the matrices~$A \in {\rm co}(\cA)$ possesses a nonnegative eigenvalue, then it is not 
Hurwitz, hence, the system is  unstable. 

\begin{figure}[h!]
\centering
 	{\includegraphics[scale = 0.77]{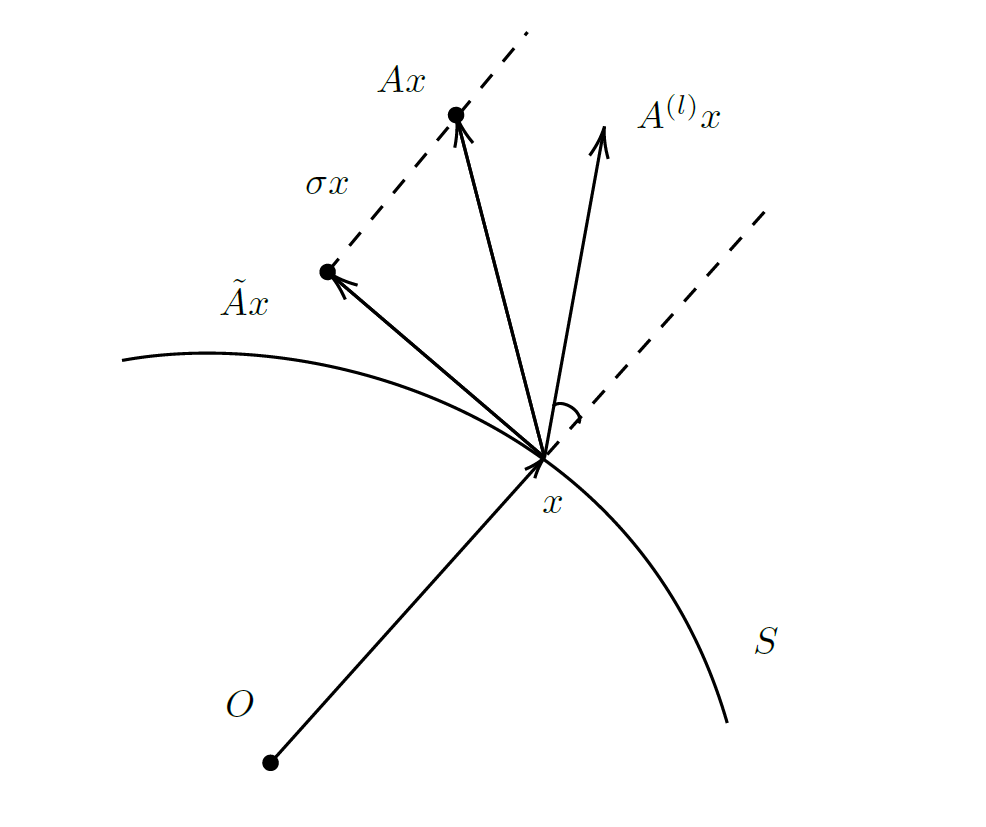}}
 	\caption{The extremal trajectory~$S$ and the vectors~$A\bx, \, \tilde A\bx, \, A^{(l)}\, \bx$}
 	\label{fig.vectors}
 \end{figure}

{\em (Sufficiency)}. Assume the contrary:   the system is unstable. If~$\sigma = 0$, 
then Theorem~\ref{th.20} implies that at least one of the leading trajectories
is well-defined and periodic, which contradicts  the assumption of the theorem. If~$\sigma > 0$,
then we consider the shifted family~$\tilde \cA = \cA - \sigma I$. 
It does not contain degenerate matrices, otherwise~${\rm co} (\cA)$ contains a
matrix with the eigenvalue~$\sigma$. Hence, by Theorem~\ref{th.20}, the system~$\tilde \cA$ 
possesses a periodic leading 
trajectory~$\tilde \bx (\vardot)$ (say, left trajectory), which goes along the  
unit sphere~$S$ of the Barabanov norm~$f$. Let~$\tilde A(\vardot)$ be the switching law 
of this trajectory. Take arbitrary~$t\ge 0$ 
and the corresponding point of the trajectory~$\tilde \bx (t) \in S$.  
Let us denote it for the sake of simplicity~$\bx = \tilde \bx(t)$, Fig.~\ref{fig.vectors}.
Let~$\alpha_{\bx}$ be the left tangent ray to the curve~$S$ at the point~$\bx$
(see definition in the Introduction). 
The vector~$\tilde A(t) \bx$
is directed along~$\alpha_{\bx}$. 
Furthermore, $\tilde A(t) \bx \, = \, A(t)\bx - \sigma \bx$, therefore, 
the vector~$A(t)\bx \, = \, \tilde A(t)  \bx\, + \, \sigma \bx$ 
lies inside the angle between the vectors~$\bx$ and $\tilde A(t) \bx$. 
Hence, the vector $A(t)\bx$ forms a smaller angle with~$\bx$ than the vector~$\tilde A(t)\bx$.
On the other hand, $A^{(l)}\bx$ forms, by definition of~$A^{(l)}$,  a still smaller angle. 
Thus, the angle between~$\bx$  and $A^{(l)}\bx$ is smaller than 
the angle between~$\bx$ and  the left tangent ray~$\alpha_{\bx}$. This means 
that the vector $A^{(l)}\bx$ starting at the point~$\bx$ is directed 
outside~$S$. This is true for every point~$\bx \in S$. Therefore, 
{\em the Barabanov norm $f(\bx(t))$ increases along the 
 left leading trajectory~$\bx(t)$}. Consequently, 
 if the system~$\cA$ has a left round leading trajectory with period~$T$ starting 
 at~$\bx(0) = \bx_0$,
 then~$f(\bx(T)) > f(\bx(0))$ and hence, $\bx(T) = - \lambda \bx_0$ for~$\lambda > 1$, i.e., this trajectory is increasing. 

It remains to consider the case 
when the left leading trajectory~$\bx(t)$ starting at some point~$\bx(0) = \bx_0$
is not round. This means that the direction of the vector~$\bx(t)$ tends to some direction~$\bx_{\infty}$
as~$t\to +\infty$, and the angle between~$\bx_0$ and~$\bx_{\infty}$ is smaller than~$\pi$. 
Let us show that in this case $\bx_{\infty}$ is an eigenvector of some operator~$A \in {\rm co}(\cA)$ with a nonnegative eigenvalue. 
The norms of the vectors~$A^{(l)}\bx$ are bounded below by some 
constant~$c>0$ in a neighbourhood of~$\bx_{\infty}$. 
Otherwise, by the compactness of~$\cA$ there is a sequence~$t_i \to \infty$
such that~$A^{(l)}(t_i)$ converges to some matrix~$\bar A \in {\rm co}(\cA)$
and~$\bar A\bx_{\infty} \, = \, 0$, which is a contradiction. 
Denote by~$h(\bx)$
the length of the orthogonal projection of~$A^{(l)}\bx$ to the line perpendicular to~$\bx$. 
We have~$\lim_{\bx \to \bx_{\infty}} h(\bx(t)) = 0$, otherwise, 
the trajectory~$\bx(t)$ reaches the direction~$\bx_{\infty}$ within finite time. 
Hence, $h(\bx_{\infty}) = 0$, so, $\bx_{\infty}$ is an eigenvector of~$\bar A$. 
If the corresponding eigenvalue is negative, then  
the vectors~$\bx_{\infty}$ and~$\bar A\, \bx_{\infty}$ have opposite directions, then 
the angle between~$\bx(t_i)$ and~$A^{(l)}\bx(t_i)$ tends to~$\pi$. 
Hence, so does the angle between~$\bx(t_i)$ and~$\tilde A(t)\bx(t_i)$, since, as we have shown above,  this angle is bigger than the angle between~$\bx(t_i)$ and~$ A^{(l)}\bx(t_i)$. 
Thus, the angle between the vector~$\bx(t_i)$ and the left tangent ray~$\alpha_{\bx(t_i)}$
 tends to~$\pi$ as~$\bx(t_i) \to \bx_{\infty}$. 
This is impossible since the angle between the vector~$\bx$ and the left tangent 
to~$S$ at the point~$\bx$ is bounded above  over all~$\bx \in S$ 
by~$\pi - \varepsilon$, where the  constant~$\varepsilon > 0$ depends only on~$S$. 
The contradiction completes the proof.

  {\hfill $\Box$}
\medskip 

Thus, if we can check that none of the operators from~${\rm co}(\cA)$ has a nonnegative 
eigenvalue, then the stability decision is reduced to construction 
of the leading trajectory. To this end, we realize the following: 
\bigskip 

\textbf{Algorithm for deciding stability}. We are given a 
compact irreducible family of $2\times 2$ matrices~$\cA$. 
The problem is to decide whether~$\sigma(\cA) < 0$.  
\smallskip 

Check whether there are matrices from~${\rm co}(\cA)$ with 
a nonnegative eigenvalue. If there are, then the system is unstable. Otherwise, 
do the following steps. 
\smallskip 

Take arbitrary~$\bx_0 \ne 0$ and find~$A^{(l)}_{\bx_0}, A^{(r)}_{\bx_0}$. 
If~$A^{(l)}_{\bx_0}$ exists, then do all the steps below. 
Then do the same for~$A^{(r)}_{\bx_0}$ provided that it exists. 
\smallskip 

Solve the ODE~$\, \dot \bx \, = \, A^{(l)}\, \bx, \ \bx(0) = \bx_0$.   Three cases are possible. 
\medskip 

1) The trajectory comes to the opposite direction, i.e.,~$\bx(t) \, = \, - \, \mu \, \bx_0$
for some~$t>0$.  If, $\mu \ge 1$, then $\sigma \ge 0$ and~$\cA$ is unstable, otherwise, the trajectory decreases. 

\smallskip 

2) The trajectory stops at some point~$\bx(\bar t), \, \bar t> 0$, 
such that the angle between the vectors~$\bx(t)$ and~$\bx (\bar t)$ is less than~$\pi$.
This means that~$A^{(l)}(\bar t)$ is not defined, i.e., no vectors~$A\bx(\bar t), \, A\in \cA$, 
form an angle in~$(0, \pi)$ with~$\bx(\bar t)$. In this case there is no left leading round trajectory.

3) The trajectory~$\bx (t)$ approaches to some direction~$\bx_{\infty}$ 
such that the angle between the vectors~$\bx(t)$ and~$\bx_{\infty}$ does not exceed~$\pi$. 
This means that the angle between the vectors~$\bx(t)$ and~$\bx_{\infty}$ tends to zero
as~$t \to +\infty$.  In this case there is no left leading round trajectory.  
\smallskip

 \smallskip 

If the case 1 does not imply instability, then do the same with the 
right leading trajectory, provided~$A^{(r)}_{\bx_0}$ exists. If the case 1 does not imply instability, then the system is stable.

 \bigskip
 
\noindent  \textbf{Comments.} By Theorem~\ref{th.30}, the system is stable if and only if~${\rm co}(\cA)$ does not contain matrix with a nonnegative eigenvalue and 
there is no round leading non-decreasing trajectory.  The latter  means that each of the two leading round trajectories (left and right) either does not exist or decreases. 
The left trajectory may not exist only if one of the cases 2 and 3 takes place. 
The same is true with the right trajectory. 

\medskip 

\begin{remark}\label{r.20}
{\em Let us emphasise that the leading trajectories of the systems~$\cA$   
and~$\tilde \cA \, = \, \cA - \sigma I$  may have different switching laws. 
Indeed, the extreme left directions 
of those systems at a given point~$\bx$ 
may correspond to different matrices from~$\cA$ (see Fig.\ref{fig.vectors}: 
the vectors~$A\bx$ and~$A^{(\ell)}\bx$ do not coincide).  This is actually the main difficulty 
in the construction of the Barabanov norm: if we find a leading trajectory of the 
family~$\cA$, it does not mean that we know the switching law 
of the trajectory of~$\tilde \cA$, which defines~$S$.} 
\end{remark}

\bigskip

\enlargethispage{2\baselineskip}

\begin{center}
\large{\textbf{4. Barabanov norms for general control sets}}
\end{center}

\medskip 

In this section we always assume that~$\sigma(\cA) = 0$. 
Theorem~\ref{th.20} proved in Section~2 gives a complete description of Barabanov norms 
when all matrices of~${\rm co}(\cA)$ are non-degenerate. In this case 
the sphere~$S$ is defined by a periodic leading trajectory. The presence of degenerate matrices
changes the situation completely  and admits a rich variety of invariant norms.   
To attack this  case, we begin with the structure of arcs of~$S$ corresponding to non-degenerate matrices. 
\medskip 


\begin{center}
\textbf{4.1. Auxiliary results}
\end{center}

\medskip

\begin{lemma}\label{l.70}
Suppose an open arc of the curve~$S$ between points~$\bx_1$ and~$\bx_2$
(enumerated counterclockwise)
does not intersect kernels of operators from~${\rm co}(\cA)$. 
 Then this arc
 is a union of two leading trajectories starting at some point~$\bs \in S$
and going to~$\bx_1$ and $\bx_2$, respectively. One of those trajectories can be empty, in which case~$\bs$ coincides with~$\bx_1$ or~$\bx_2$. 
\end{lemma}
{\tt Proof.} 
Moving the points~$\bx_1, \bx_2$ slightly inside the arc  
we obtain an embedded closed arc~$\bx_1'\bx_2'$. If we prove the statement 
for this arc, then it will follow for the original open arc by the limit passage. 
The remainder of the proof is the same as the proof Theorem~\ref{th.20}. 
First, we show that for every~$\bx$ on the arc~$\bx_1'\bx_2'$, the angle between the 
leading rays is less than~$\pi$. If this is not the case, then there exists 
a convex combination of the operators~$A^{(l)}_{\bx}, A^{(r)}_{\bx}$
for which~$\bx$ is an eigenvector with a nonnegaive eigenvalue~$\lambda$. 
If~$\lambda > 0$, then~$\sigma(\cA) > 0$. If~$\lambda = 0$, 
then~$\bx$ belongs to the kernel of that convex combination, which contradicts the assumptions. 
Then we prove that  the arc~$\bx_1'\bx_2'$
contains only a finite number of switching points between the leading directions.
If moving along the arc~$\bx_1'\bx_2'$ counterclockwise 
we change the leading direction from~$a_{\bx}$ to~$b_{\bx}$, i.e, 
we meet an arc~$\bz_1\bz_2$ with the leading direction~$a_{\bx}$
followed by an arc~$\bz_2\bz_3$ with the leading direction~$b_{\bx}$, 
then as in the proof of Step 3 of Theorem~\ref{th.20} we come to the 
contradiction. Thus, all the switching points 
change the leading direction from~$b_{\bx}$ to~$a_{\bx}$. 
Hence, either there are no switching points or there is a unique 
switching point~$\bs$ from which two trajectories start 
and go to the ends of the arc.

  {\hfill $\Box$}
\medskip 

The point~$\bs$ will be referred to as a {\em source} for the arc~$\bx_1\bx_2$. 
Lemma~\ref{l.70} asserts that every arc which does not have interior points on the kernels
of operators from~${\rm co}(\cA)$ contains a unique source, possibly coinciding with an end of the arc. For constructing an invariant norm, we need an inverse problem: 
under which conditions does the arc~$\bx_1\bx_2$ composed of two leading trajectories exist? 
More precisely, let~$\ell_1, \ell_2$ be rays starting at the origin and the angle between them 
does not exceed~$\pi$ (the angle is oriented, the enumeration is counterclockwise). 
For given points~$\bx_i \in \ell_i \setminus \{0\}, \, i=1,2$, 
one needs to determine whether there exists a source point~$\bs$ and two leading trajectories
going to~$\bx_1$ and~$\bx_2$, respectively? 

To formulate the criterion, we need some further notation.  Assume the left leading trajectory 
goes from a point~$\bz_1 \in \ell_1$ to~$\bx_2 \in \ell_2$. Denote~$k_1 = \frac{\|\bx_2\|}{\|\bz_1\|}$. The right leading trajectory 
goes from~$\bz_2 \in \ell_2$ to~$\bx_1 \in \ell_1$. Denote~$k_2 = \frac{\|\bx_1\|}{\|\bz_2\|}$.
Finally, let~$m\, = \, \min\, \bigl\{\ln k_1, - \ln k_2\bigr\}$ and~$M
= \, \max\, \bigl\{\ln k_1, - \ln k_2\bigr\}$. 
\begin{lemma}\label{l.75}
Let~$\sigma(\cA) = 0$, 
$\ell_1, \ell_2$ be rays starting at the origin, the  angle between them 
does not exceed~$\pi$ and 
does not contain  the kernels of operators from~${\rm co}(\cA)$
different from~$\ell_1, \ell_2$. 
Suppose there exists both (left and right) leading trajectories connecting points 
of those rays; then for given~$\bx_i \in \ell_i, \, i=1,2$, there are two trajectories with a common source leading to~$\bx_1$ and~$\bx_2$, respectively, if and only if 
\begin{equation}\label{eq.mM}
m \quad \le \quad  \ln \ \frac{\|\bx_2\|}{\|\bx_1\|} \quad  \le \quad  M, 
\end{equation}
where~$m, M$ are defined above. 
\end{lemma}
{\tt Proof.}  The source point~$\bs$ exists if and only if 
the leading trajectories~$\bz_1\bx_2$ (from~$\ell_1$ to~$\ell_2$) 
and~$\bz_2\bx_1$ (from~$\ell_2$ to~$\ell_1$)  intersect.
This happens when either~$\|\bx_1\| \le \|\bz_1\|, \, \|\bx_2\| \le \|\bz_2\|$
or conversely:~$\|\bx_1\| \ge \|\bz_1\|, \, \|\bx_2\| \ge \|\bz_2\|$. 
In the former case we have~$\frac{\|\bx_2\|}{\|\bx_1\|} \, \ge \, 
\frac{\|\bx_2\|}{\|\bz_1\|}\, = \, k_1$ and 
~$\frac{\|\bx_2\|}{\|\bx_1\|} \, \le \, 
\frac{\|\bz_2\|}{\|\bx_1\|}\, = \, \frac{1}{k_2}$. Thus, $\frac{\|\bx_2\|}{\|\bx_1\|} \in 
\bigl[k_1 , \frac{1}{k_2} \bigr]$. This is possible when~$k_1 \le \frac{1}{k_2}$.  Similarly, the latter case means that
$\frac{\|\bx_2\|}{\|\bx_1\|} \in 
\bigl[\frac{1}{k_2}, k_1\bigr]$, which is equivalent to~(\ref{eq.mM}).

  {\hfill $\Box$}
\medskip

\begin{lemma}\label{l.65}
If~$\sigma(\cA) = 0$ and there is a degenerate matrix~$A \in {\rm co}(\cA)$, 
then its second eigenvalue is negative. 
\end{lemma}
{\tt Proof.}  The second eigenvalue cannot be positive; if it is zero, then 
either~$A=0$, this case is excluded, or the equation~$\dot \bx = A\bx, \, \bx(0) = \bx_0$, 
  has a solution~$\bx(t) = \bx_0 + t A\bx_0$, which is unbounded whenever~$\bx_0 \notin {\rm Ker}\, A$. This contradicts the existence of the invariant norm.

  {\hfill $\Box$}
\medskip

Let us remember that~$\alpha_{\bx}$ and~$\beta_{\bx}$ denote, respectively,
the left and right tangent rays to~$S$ at the point~$\bx\in S$.

\begin{prop}\label{p.60}
Suppose an operator~$A \in {\rm co}\, (\cA)$ is degenerate and its kernel meets~$S$
at a point~$\bx_0$; then  
 the line passing through~$\bx_0$ parallel to~${\rm Im}\, A$ is a line of support 
for~$S$. 
\end{prop}
{\tt Proof.}  Denote by~$\ell$ the line passing through~$\bx_0$ parallel to~${\rm Im}\, A$. If~$\ell$ is not a line of support, then there is a point~$\by \in S$
 separated from the origin by~$\ell$. Draw a line~$\ell'$ parallel to~$\ell$ through~$\by$ and 
denote by~$\by'$ the point of its intersection with the line~$O\bx_0 \, = \, {\rm Ker}\, A$. 
Clearly, $|O\by'|\, > \, |O\bx_0|$, hence, $\by'$ lies outside~$S$. 
The trajectory~$\bx(t)$ of the equation~$\dot \bx\ = \, A\bx, \, \bx(0) = \by,$ approaches~$\by'$
along 
the line~$\ell'$ as~$t\to +\infty$. 
Therefore,~$\bx(t)$ does not stay  inside~$S$, which contradicts 
 the definition of Barabanov's norm.

  {\hfill $\Box$}
\medskip

\begin{lemma}\label{l.80}
Suppose under the assumptions of Lemma~\ref{l.70}, 
we have~$\bs \ne \bx_1$ and~$\bx_1 \in {\rm Ker}\, A$ for some $A\in {\rm co}(\cA)$; 
then the ray~$\beta_{\bx} = b_{\bx}$, which defines the trajectory from~$\bs$ to~$\bx_1$, 
tends to~$\, - \, \alpha_{\bx_1}$ as~$\bx \to \bx_1$. Moreover, there exists~$A'\in {\rm co}(\cA)$
with the kernel containing~$\bx_1$ and image parallel to~$\alpha_{\bx_1}$. 
\end{lemma}
{\tt Proof.} Since the arc~$\bx_1\bs$ does not contain switching points, it belongs to~$C^1$. 
Hence, at every point~$\bx$ of this arc, the tangent line is well-defined and 
contains~$\beta_{\bx}$. By the convexity of~$S$, this tangent line 
converges to the left tangent at~$\bx_1$, i.e., to~$\alpha_{\bx_1}$ when $\bx \to \bx_1$. 
This proves that~$\beta_{\bx} \, \to \, - \alpha_{\bx_1}$. 
 Let~$\ell$ be the line passing through~$\bx_1$ parallel to the image of~$A$. By Proposition~\ref{p.60}, $\ell$ is a line of support for~$S$. 
If $\alpha_{\bx_1}$ lies on~$\ell$, then everything is proved. Otherwise, $\ell$
separates the ray~$-\alpha_{\bx_1}$ from~$S$. Hence, $-\alpha_{\bx}$ is directed outside~$S$. 
Furthermore, the ray~$b_{\bx_1} = \lim_{\bx \to \bx_1}b_{\bx}$ 
has the same direction as~$-\alpha_{\bx_1}$. 
Hence, if~$A^{(r)}\bx_1 \ne  0$, then this vector  goes outside~$S$, 
which is impossible. 
If  $A^{(r)}\bx_1 \, = \, 0$, then an arbitrary 
partial limit~$A' = \lim_{\bx_i \to \bx_1} A^{(r)}_{\bx_i}$ 
satisfies~$A'\bx_1 = 0$ and its image 
contains the limit direction~$\lim_{\bx \to \bx_1} b_{\bx} \, = \, -\alpha_{\bx_1}$. 

  {\hfill $\Box$}
\medskip

Let~$A_1, A_2$ be two degenerate operators with negative second eigenvalues and  with a common image~$L \subset \re^2$ but with different kernels~$K_1, K_2$. 
The convex hull~$\bigl\{(1-\lambda)A_1 + \lambda A_2: \ \lambda \in [0,1] \bigr\}$
is called a {\em reverse pencil}~$\cP[A_1, A_2]$. 
A segment~$[\bx_1, \bx_2]$ with ends on~$K_1,K_2$ and parallel to~$L$ is called a~{\em reverse segment}. 
Then for every~$\bx \in [\bx_1, \bx_2]$, the images~$A_1\bx, A_2\bx$ are 
parallel to~$L$ and have opposite directions. 
\begin{lemma}\label{l.35}
If a family~$\cA$ with~$\sigma(\cA) = 0$ contains two 
degenerate operators~$A_1, A_2$ with a common image, 
then the pencil~$\cP[A_1, A_2]$ lies in~${\rm co}(\cA)$ and  
$S$ contains a reverse segment with ends on the kernels. 
\end{lemma}
{\tt Proof.} Using the notation above, we consider the set of reverse segments 
intersecting~$S$ and denote by~$[\bx_1, \bx_2]$ the one most distant from the origin. 
Denote also by~$\by_i$ the intersection of~$S$ with the ray~$O\bx_i, \, i = 1,2$. 
The arc~$\by_1\by_2 \subset S$ lies inside the triangle~$\bx_1O\bx_2$ 
and has a common point~$\bx_0$ with the segment~$[\bx_1, \bx_2]$. 
The trajectories of the 
systems~$\dot \bx = A_1\bx$ and~$\dot \bx = A_2\bx$ starting 
at~$\bx_0$ fill the half-intervals
~$[\bx_0, \bx_1)$ and~$[\bx_0, \bx_2)$,  respectively.  Both those trajectories are inside~$S$ 
by the Barabanov norm property. Hence, the whole reverse segment~$[\bx_1, \bx_2]$ is inside~$S$. 
Therefore, it coincides with the arc~$\by_1\by_2$. 

  {\hfill $\Box$}
\medskip

Our main auxiliary result states that the leading trajectory on~$S$ may change its direction only on 
the reverse intervals unless~$\cA$ has  complex dominance (the ellipsoidal case).  
  
\begin{lemma}\label{l.40}
If there is a point~$\bx$  for which~$a_{\bx}$ and~$b_{\bx}$ are collinear and have opposite directions, then $\sigma \ge 0$. If $\sigma = 0$, then either $\cA$ has complex 
dominance or~${\rm co}(\cA)$ contains a reverse pencil. In the latter case 
$\bx$ is an interior point of the reverse segment on~$S$. 
\end{lemma}
The proof is in Appendix. 
\smallskip

If~$\sigma(\cA) =0$ and ${\rm co}(\cA)$ contains a degenerate 
operator, whether so does~$\cA$? In general,  the answer is negative: 
if a matrix~$A$ has an imaginary spectrum, i.e., defines an ellipsoidal rotation, 
then so does~$-A$ and the family~$\cA = \{A, -A\}$
satisfies~$\sigma(\cA) = 0$ and  
  ${\rm co}(\cA)$ contains the zero matrix. It turns out, however, that apart from this case, 
  the answer is affirmative.  As usual, we call a point of a convex set {\em extremal}
  if it is not a midpoint of a segment in this set. 

\begin{theorem}\label{th.35}
Suppose~${\rm co}(\cA)$ contains a dominant matrix 
(with real or complex dominance); then so does the set of extreme points of~${\rm co}(\cA)$. 
\end{theorem}
The theorem immediately follows from a more general fact below: 
\begin{prop}\label{p.50}
If~$\sigma(\cA) =0$ and ${\rm co}(\cA)$ contains a degenerate matrix~$A$, 
which is not in~$\cA$, then  one of the three following conditions are satisfied: 
\smallskip 

1) $\cA$  has  complex dominance;  
\smallskip 

2) there are operators $A_1, A_2 \in \cA$  such that~$A \in \cP[A_1,A_2]$; 
\smallskip 

3)  $A$ is a linear combination of two degenerate operators from~$\cA$
with the same kernel. 
\end{prop}
The proof is in Appendix. 
\begin{cor}\label{c.55}
Assume that~$\sigma(\cA) =0$ and $\cA$~does not have complex dominance. 
Let a matrix~$A \in {\rm co}(\cA)$ be degenerate with a kernel~$K$. 
Then either~$A$ belongs to a reverse pencil or there is a matrix from~$\cA$
with the same kernel~$K$. 
\end{cor}
{\tt Proof}. Proposition~\ref{p.50} yields that either $A\in \cA$
(in which case the proof is completed) or $A$ satisfies one of the cases 
2), 3) of that proposition. The latter means that  either $A \in \cP\, [A_1, A_2]$ for some~$A_1, A_2 \in \cA$ or~$A$ is a linear combination of two 
 operators from~$\cA$ with kernel~$K$.

  {\hfill $\Box$}
\medskip

\newpage

\begin{center}
\textbf{4.2. The structure of the Barabanov norms}
\end{center}

\medskip

\begin{figure}[h!]
\centering
 	{\includegraphics[scale = 0.9]{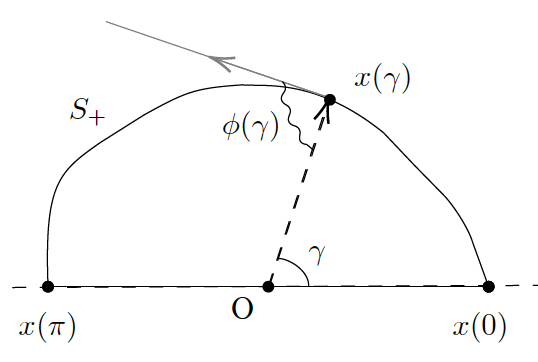}}
 	\caption{Parametrization of the curve~$S_+$}
 	\label{fig.curve}
 \end{figure}

Now we can give a complete classification of all Barabanov norms
for a set of~$2\times 2$ matrices~$\cA$. 
As usual we assume that~$\cA$ is irreducible and~$\sigma(\cA) = 0$. 
The unit sphere~$S$ of an invariant norm~$f$ is parametrized 
by the angle~$\gamma$ (the {\em argument} of a point~$\bx\in S$) 
 between the vector~$\bx$ and~the~$OX$ axis. 
 Due to the symmetry, 
we consider only the upper half-sphere~$S_{+}\, = \, 
\{\bx(\gamma): \ \gamma \in [0,\pi]\}$, where~$\bx(\pi)= - \bx(0)$. 
The angle between the vector~$\bx(\gamma)$ and the left tangent to~$S$
at the point~$\bx$ is denoted as~$\varphi(\gamma)$, Fig.~\ref{fig.curve}. 

We consider the function~$g(\gamma) = 
\ln \, \|\bx(\gamma)\|$, where the norm is Euclidean. Clearly,~$g(\cdot)$ is Lipschitz continuous 
on~$[0, \pi]$ and it defines~$S$ in a unique way. 
\begin{lemma}\label{l.60}
The derivative~$g'(\gamma)$ exists for almost all~$\gamma$ (in the Lebesgue measure on~$[0,\pi]$) and is equal 
to~$\cot \varphi(\gamma)$. 
\end{lemma}
{\tt Proof.} Since the set of corner points of~$S$ is at most countable, 
it has a tangent at~$\bx(\gamma)$ for almost all~$\gamma$.  
For a small~$\delta$, consider the triangle with vertices~$O, \bx(\gamma), 
\bx(\gamma + \delta)$. Applying the sine law and 
taking into account that the angles at the 
vertices~$\bx(\gamma), \bx(\gamma + \delta)$ are close to~$\pi - \varphi(\gamma)$ and~$\varphi(\gamma) - \delta$,
respectively, we get
$$
g(\gamma +\delta)\, - \, g(\gamma) \ = \ 
\ln\, \frac{\|\bx(\gamma + \delta)\|}{\|\bx(\gamma)\|} \ = \ 
\ln \, \frac{\sin \varphi(\gamma)}{\sin (\varphi(\gamma) - \delta)}  \ + \ o(\delta) \ = \ 
\delta \, \cot \varphi(\gamma) \  + \  o(\delta) \quad \mbox{as}\quad \delta \ \to \ +0, 
$$
from which the lemma follows. 

  {\hfill $\Box$}
\medskip 

Note that if~$\bx$ belongs to the kernel of some operator~$A\in {\rm co}\,(\cA)$, 
then $\bx + {\rm Im}\,(A )$ is a line of support of~$S$ (Proposition~\ref{p.60}). 
Since $S$ has at most a countable set of corner points and at all other points the 
line of support is unique, we obtain:  
\begin{prop}\label{p.80}
Denote by~$S'$ the set of points~$\bx \in S_+$
such that~$\bx$ belongs to the kernel of some~$A_{\bx}\in {\rm co}\,(\cA)$. 
Then for almost all points of~$S'$ the angle between~$\bx$ and~${\rm Im}\,(A_{\bx})$, 
 where~$\bx = \bx(\gamma)$, is equal to~$\varphi(\gamma)$. 
\end{prop}
\medskip

The curve~$S_+$ is a disjoint union of three sets (in what follows we identify the number~$\gamma \in [0,\pi]$ and the point
$\bx(\gamma) \in S_+$): 
\smallskip 

1) The reverse open intervals. We denote by~$R$ the set of those intervals. 
\smallskip 

2) The set~$D$ of points   which belong to 
the kernels of  operators from~$\cA$ and do not lie on the 
 reverse open intervals.  
For each~$\bx(\gamma) \in D$, the line~$\{\lambda \bx: \ \lambda \in \re\}$
is the kernel of some~$A \in \cA$. By Proposition~\ref{p.60}, $\bx + {\rm Im} (A)$ is 
a line of support for~$S$. By Corollary~\ref{c.55}, every point~$\bx(\gamma) \notin D$, 
does not belong to kernels of operators from~${\rm co}(\cA)$.  
\smallskip 

3) The set of points, which are neither from reverse intervals nor from~$D$. 
This set is open in~$S_+$ (as a complement to the closed set) and hence,  is a disjoint union of 
open intervals.  
The set of those intervals is countable. 
 Each interval from this set 
does not intersect kernels of operators from~$\cA$ and therefore, 
 admits two leading trajectories (Lemma~\ref{l.70}).    The set of intervals with a unique leading trajectory (when the second trajectory is empty and the source point~$\bs$ coincides with one of the ends of the interval) is denoted by~$H$. The set of intervals 
 with two leading trajectories is denoted by~$P$. 
 Thus, $P$ and~$H$ are disjoint sets of intervals, each interval does not intersect kernels of 
 operators from~$\cA$. 
 
Let~$Q = R\cup H$. The sets~$P$ and $Q$ are finite or countable, we enumerate their 
 elements as~$p_i$ and~$q_i$ respectively. 
Clearly, those  sets are not ordered on the segment~$[0,\pi]$). 

\begin{figure}[h!]
\centering
 	{\includegraphics[scale = 0.73]{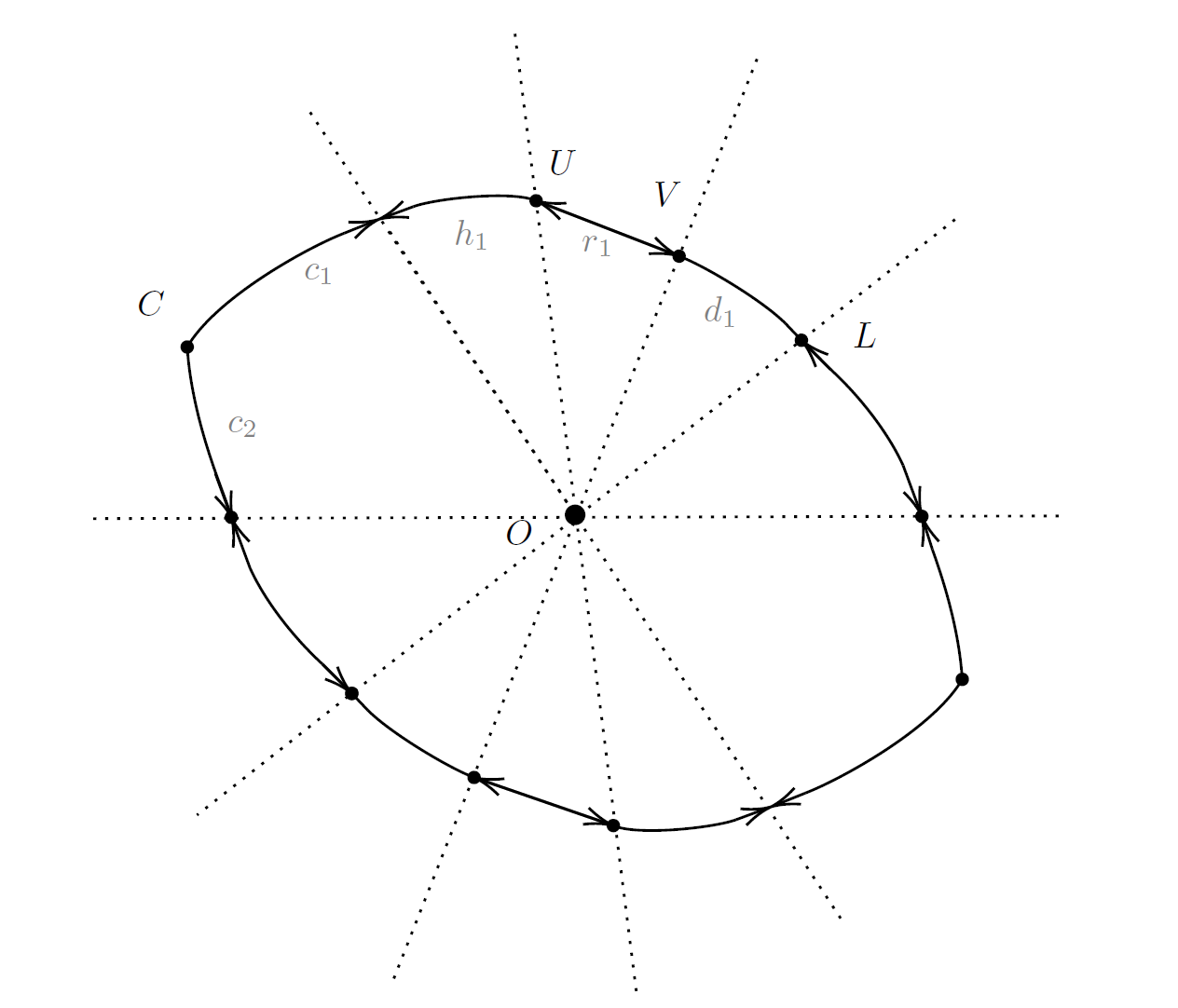}}
 	\caption{The  invariant sphere~$S$}
 	\label{fig.classification}
 \end{figure}

Fig.~\ref{fig.classification}
shows different types of arcs of the curve~$S$. The reverse interval~$r_1$ belongs to~$R$;  
the operators of the corresponding pencil have a common image~$UV$, their kernels 
run from the line~$OU$ to~$OV$.   
The  arcs $c_1$ and~$c_2$ are the leading trajectories going from a common source~$C$, 
they form an interval~$p_1\in P$. The arc~$h_1$ belongs to~$H$, 
it is a leading trajectory going  from the source~$U$. 
The arc~$d$ is a part of the set~$D$, all points of this arc belong to kernels of operators 
from~$\cA$. 
\medskip 

Without loss of generality, with possible change of the coordinates,
  it can be assumed that 
$\bx(0)$ does not belong to the open set~$Q \cup P$. 

For every interval~$q=(\gamma_1, \gamma_2)$, we denote~$\Delta(q) = g(\gamma_2) - g(\gamma_1)$. 
For each~$q_i \in Q$, the value~$\Delta(q_i)$ is well-defined and easily found; 
  for~$p_j \in P$, the value~$\Delta(p_j)$ is arbitrary 
  from the segment~$[m_j, M_j]$, where the values~$m_j, M_j$
  are defined in Lemma~\ref{l.75} under the names~$m, M$. 
  Inequality~(\ref{eq.mM}) implies that~$\Delta(p_j) \in [m_j, M_j]$. 

\begin{theorem}\label{th.40}
Let~$\cA$ be a control set of~$2\times 2$ matrices of system~(\ref{eq.main}) 
with~$\sigma(\cA) = 0$
and $S$ be the unit sphere of its Barabanov norm. 
Then  for every~$\gamma \in [0,\pi]$ which is either an element of~$D$ or an end of an interval 
from~$Q$ or~$P$, we have 
\begin{equation}\label{eq.g}
g(\gamma)\ = \ g(0) \ + \ \int_{[0, \gamma]\cap D} \cot \varphi(\tau)\, d\tau \ + \ 
\sum_{q_i \in Q, \, q_i \subset [0, \gamma]}\, \Delta(q_i) \ + \ 
\sum_{p_j \in P, \, p_j \subset [0, \gamma]}\, \Delta(p_j)\, . 
\end{equation}
For given numbers~$s_j \in [m_j, M_j]$, there exists 
an invariant sphere~$S$ with~$\Delta(p_j) = s_j, \, p_j \in P$, 
if and only if 
\begin{equation}\label{eq.si}
\sum_{p_j \in P}\, s_j  \ = \ 
-  \int_{D} \cot \varphi(\tau)\, d\tau \ - \ 
\sum_{q_i \in Q}\, \Delta(q_i) \, . 
\end{equation}
where~$\varphi(\tau)$ is the angle between the vector~$\bx(\tau)$ 
and the image of the corresponding matrix~$A$ whose kernel contains~$\bx$. 
For every such sequence~$\{s_i\}$, the invariant sphere~$S$ is unique up to normalization. 
\end{theorem}
{\tt Proof.} The function~$g(\tau)$ is  Lipschitz and hence, absolutely 
continuous on~$[0,\pi]$. Therefore, it is a primitive of its derivative, which is equal to~$\cot \varphi(\tau)$ (Lemma~\ref{l.60}). 
For almost all points~$\tau \in D$, there exists a tangent line to~$S$
at the point~$\bx(\tau)$, hence, it is parallel to the image 
of the corresponding matrix~$A$~(Proposition~\ref{p.60}). 
Thus, for almost all~$\tau \in D$, the angle between~$\bx(\tau)$ and~${\rm Im} (A)$ 
is equal to~$\varphi(\gamma)$. 

Substituting~$\gamma = \pi$ in the equation~(\ref{eq.g})
and taking into account that~$g(\pi) = g(0)$ because of the symmetry of~$S$, we obtain~(\ref{eq.si}). 

Assume now that we have a set of numbers~$s_j \in [m_j, M_j]$
satisfying~(\ref{eq.si}). Define the function~$g(\gamma)$ by the 
formula~(\ref{eq.g}) for all~$\gamma \in D$ and for all~$\gamma$ that are  ends of the intervals~$p_i, q_j$. 
Let us recall that we identify each interval with the corresponding angle 
centered at~$O$. This defines the value~$\bx(\gamma)$ at all points of the set~$D$ and
all ends of the intervals~$q_i, p_i$. Then~$\bx(\gamma)$ is extended in a unique way 
to each interval~$q_i \in Q$. Indeed, if this interval (denote it by~$(\gamma_1, \gamma_2)$) belongs to~$R$, then  the arc of~$S$ in this angle is a segment connecting the 
points~$\bx(\gamma_1)$ and~$\bx(\gamma_2)$; if this interval 
belongs to~$H$, then  the arc of~$S$  is a 
unique leading trajectory connecting those points. Finally, on 
each  interval~$p_i$ (denote it again  by~$(\gamma_1, \gamma_2)$), the arc of~$S$ is a union of two arcs of the leading trajectories starting at~$\bx(\gamma_1)$ and~$\bx(\gamma_2)$. 
Those trajectories intersect since~$g(\gamma_2) - g(\gamma_1) = s_i \in [m_i, M_i]$ 
(Lemma~\ref{l.75}).

  {\hfill $\Box$}
\medskip 

\medskip

\bigskip 

\begin{center}
\textbf{4.3. Construction and the  uniqueness}
\end{center}

\medskip

Theorem~\ref{th.40} is constructive, it allows us to 
obtain the Barabanov norm algorithmically. 
\smallskip 

\textbf{Algorithm.}
To construct a norm, one needs first 
to find~$\sigma(\cA)$.  This can be done by bisection with 
the comparison to zero (i.e. deciding the stability) in each step. 
Then we shift the family so that~$\sigma(\cA) = 0$. 
Find all degenerate operators in~${\rm co}\, (\cA)$. 
If there are none, then one of the 
leading trajectories (left or right) of~$\cA$ is periodic and its period 
describes the invariant sphere~$S$ (Theorem~\ref{th.20}). 
Otherwise, we find the set of reverse pencils (it corresponds to the set~$R$ of reverse intervals).  This gives us the directions of all reverse segments on~$S$. 
 Thus, we know~$\Delta(r_i)$ for all~$r_i \in R$. 
 Then we find all intervals that are not intersected by 
 kernels of operators from~$\cA$ and construct the leading trajectories for each of them. 
This construction is unique up to homothety with respect to the origin~$O$.  
On each of those  intervals~$(\gamma_1, \gamma_2)$ the leading  trajectories are uniquely defined by the values~$g(\gamma_i) \, = \, \ln \|\bx(\gamma_i)\|, \, i=0,1$, which will be found later. 
The intervals  having one leading trajectory are denoted by~$h_j  \subset [0,\pi]$ and 
others are denoted by~$p_i \subset [0,\pi]$.  
Compute~$\Delta(h_j)$. This value does not depend on the 
coefficient of homothety, i.e., on the value~$\|\bx(\gamma_1)\|$ at the end of this interval. 
We obtain the 
term~$\sum\limits_{q_i \in Q, \, q_i \subset [0, \gamma]}\, \Delta(q_i)$  
in the sum~(\ref{eq.g}). Then we choose an arbitrary sequence~$\{s_j\}$ satisfying~(\ref{eq.si}) and such that~$s_i\in [m_i, M_i]$ for all~$i$. If such a sequence does not exist, then by Theorem~\ref{th.40}, there is no Barabanov norm, which is a contradiction.  
Thus, choosing a sequence~$s_i$ we  substitute~$\Delta(p_i) = s_i$ into the formula~(\ref{eq.g}). To apply this formula it remains 
to define~$\varphi(\gamma)$ for~$\gamma \in D$. 
For all degenerate operators from~$\cA$ whose kernels
do not intersect the intervals from~$Q$ and~$P$,  we find the angle~$\varphi$
between the kernel and the image. The set~$D = [0,\pi]\setminus (Q\cup P)$ is filled by the kernels of those operators and 
parametrized by the angle~$\gamma$. For almost all~$\gamma \in D$, the angle~$\varphi(\gamma)$  
is unique, hence the integral in~(\ref{eq.g}) is well-defined and computed for 
all admissible~$\gamma$. Thus, we find~$g(\gamma) = \ln \|\bx(\gamma)\|$
for all~$\gamma$ from~$D$ and at the ends of the intervals~$p_i, q_i$. 
Then, having obtained~$\bx(\gamma)$ at the ends of these intervals, we extend the 
function~$\bx(\cdot)$ inside them  in a unique way. This completes the construction of~$S$.

\medskip 

The construction implies that every set of numbers~$\{s_j \in [m_j, M_j]\}$
satisfying~(\ref{eq.si})
generates a unique Barabanov norm. Hence,~$\cA$ possesses a unique Barabanov norm 
if this sequence is unique.

\begin{theorem}\label{th.50}
The Barabanov norm of second-order system~(\ref{eq.main}) is unique 
if and only if one of the following conditions is satisfied: 

1) $\cA$ has at most one degenerate operator; 

2) the set~$P$ has at most one element; 

3) the right hand side of equality~(\ref{eq.si})
is equal either to $\sum_j M_j$ or to $\sum_j m_j$. 

Otherwise, there are infinitely many Barabanov norms.  
\end{theorem}
{\tt Proof.} The uniqueness of the Barabanov norm is equivalent to the uniqueness of the 
feasible sequence~$\{s_j\}$. It is certainly unique in cases~2) and~3). The case 1) is simply reduced to the case 2) by the choice of coordinates. If~$\cA$ does not contain degenerate operators, then 
the uniqueness follows by Theorem~\ref{th.20}; if it has a unique 
degenerate operator, then we choose the axis $OX$ passing through 
its kernel, after which the set~$P$ becomes one-element (the only interval is the entire~$[0,\pi]$).  Conversely, if the sum~$\sum_j s_j$ contains at least two terms 
and it is strictly between~$\sum_j M_j$ and~$\sum_j m_j$, then there 
exist two indices (call them $1$ and~$2$) such that
~$s_1 < M_1$ and~$s_2 > m_2$. Replacing them by~$s_1 + \varepsilon$ 
and~$s_2 - \varepsilon$ we obtain a new feasible sequence, whenever~$\varepsilon > 0$
is small enough.

  {\hfill $\Box$}
\medskip

\newpage

\begin{center}
\large{\textbf{5. Each norm in~$\re^d$ is Barabanov  for a suitable system}}
\end{center}
\bigskip

The following result will be established not only for $d=2$ but in arbitrary dimension~$d$.

\begin{theorem}\label{th.60}
Every norm in~$\re^d$ is a Barabanov norm for a suitable linear switching system. 
If the unit ball of this norm is a polyhedron with~$2n$ facets, then there is the 
corresponding system with an~$n$-element control set.   
\end{theorem}
{\tt Proof.} Let~$G$ be the unit ball of this norm and~$S=\partial G$. 
For an arbitrary~$\bx \in S$, we denote by~$V_{\bx}$ the hyperplane parallel to the plane of support of~$S$
at the point~$\bx$. If this plane is not unique we choose one of them. 
Then the operator~$A_{\bx}$ is defined by the equalities~$A_{\bx}(\bx) = 0, \, 
A_{\bx}(\by) = - \by$ for all~$\by \in V_{\bx}$. 
Let us show that the family~$\cA = \{A_{\bx}: \ \bx \in S\}$
has the Barabanov norm with the unit sphere~$S$.  Indeed, for~$\bz = \bx$, the vector 
 $0 = A_{\bx}\bz$ starting at~$\bz$ is tangent to~$S$ and for every other~$\bz \in S$, the vector~$A_{\bx}\bz$ starting at~$\bz$
 is either tangent to~$S$ or directed inside~$S$. To show this we write~$\bz = \gamma \bx + \by$ for some 
 $\gamma \in [-1, 1],  \, \by \in V_{\bx}$. Then, for a small~$t> 0$, we have 
 $\bz + t A_{\bx}\bz \, = \, \gamma \bx + \by  - t \by\, = \, (1-t)\bz + t \gamma \bx
\in \, G $ since both~$\bz$ and~$\gamma \bx$ belong to~$G$.  

If $G$ is a polyhedron, then we enumerate the pairs of 
its opposite facets by~$1, \ldots , n$ and to the~$i$th pair 
 associate  an interior  point~$\bx_i$ of one of those two facets. 
 Then the~$n$ operators~$A_{\bx_i}, \, i=1, \ldots, n$, defined as above 
 have~$G$ as a unit ball of the invariant norm.

  {\hfill $\Box$}
\medskip

Now we turn back to the case~$d=2$. The system~$\cA$ constructed in  Theorem~\ref{th.60} consists of 
degenerate operators~$\cA$. This is unavoidable: 
there are norms in~$\re^2$ that cannot be invariant for families of 
non-degenerate operators with~$\sigma(\cA) = 0$.  The 
following proposition states that every non-smooth norm possesses this property. 
 
\begin{prop}\label{p.70}
Let~$G \subset \re^2$ be a unit ball of the Barabanov norm 
of a family~$\cA$  with~${\sigma(\cA) = 0}$. 
Then every arc on the sphere~$S$ connecting two corner points 
intersects an image of some degenerate operator from~${\rm co}(\cA)$. 
\end{prop}
{\tt Proof.} Consider a closed arc connecting two  corner points~$\bx_1, \bx_2$
(the enumeration is counterclockwise). 
If it does not intersect kernels of operators from~${\rm co}(\cA)$, 
then we can apply Lemma~\ref{l.70} and conclude that this arc consists of two leading trajectories going from a point~$\bs$ to~$\bx_1$ and~$\bx_2$. At least one of them, say, $\bs\bx_2$ is nontrivial. Since the arc~$\bs\bx_2$ is~$C^1$, 
it follows that at every its interior point~$\bx$, the 
ray~$a_{\bx}$
is tangent to~$S$ and is directed along~$- \beta_{\bx_2}$,  where~$\beta_{\bx_2}$ is the 
right tangent ray. 
On the other hand, $a_{\bx}$ continuously depends on~$\bx$. 
In the limit as~$\bx \to \bx_2$, we obtain: 
$a_{\bx_2}$ is directed along~$- \beta_{\bx_2}$. 
However, 
$\bx_2$ is a corner point, hence, the vector~$- \beta_{\bx_2}$  is directed outside~$G$. 
Thus,~$a_{\bx_2}$ is directed outside~$S$, which is impossible for the Barabanov sphere.

  {\hfill $\Box$}
\medskip


\begin{remark}\label{r.70}{\em 
After Proposition~\ref{p.70} a natural question arises:  
if every smooth norm is invariant for a suitable set of
{\em non-singular} matrices?  The answer is still negative. 
To see this, consider a strictly convex body~$G$ with some corner points. 
Every family~$\cA$ for which~$G$ is an invariant body includes degenerate 
operators. The polar~$G^*$, which does not have corner points,
is an invariant body for the set of adjacent operators~$\cA^*$, see~\cite{Mu2}.  
Clearly, $\cA^*$ also includes degenerate operators.  

Moreover, there are smooth strictly convex bodies 
that are not invariant bodies of families of non-singular matrices. 

}
\end{remark}

\begin{center}
\large{\textbf{6. Finite control sets and matrix polytopes}}
\end{center}
\bigskip 

If $\cA = \{A_1, \ldots , A_m\}$ is finite, then ${\rm co}\,(\cA)$ is a polytope in the space of $2\times 2$ matrices 
with vertices from~$\cA$. In this case every leading trajectory 
consists of finitely many arcs, each corresponds to one matrix. 
The {\em switching point}~$\bx$ between the arcs corresponding to matrices~$A_i$ and $A_j$
is characterized by the property that the vector~$A_i\bx$ is directed along $A_j\bx$, i.e., 
$(A_i - \lambda A_j)\bx = 0$, where~$\lambda \ge 0$ (the case~$\lambda = +\infty$ is allowed and means that~$A_j\bx=0$). Hence, to construct the leading trajectory one needs to 
solve the equation~${\rm det}\, (A_i - \lambda A_j) = 0$ for all pairs~$(i,j)$ and thus, find all 
switching points~$\bx_{ij}$, which define the whole trajectory.  We start with an arbitrary point~$\bx_0\ne 0$ and go along the left leading trajectory. Assume~$A^{(l)}_{\bx_0} = A_1$. 
 Then we go along the trajectory~$e^{\, t A_1}\bx_0$ to the closest switching point
 with the first index~$i=1$. Let it be~$\bx_{1j}$. Then we switch to the trajectory
 ~$e^{\, t A_j}\bx_{1j}$ to the closest switching point with the first index~$j$, etc. 
In total there will be at most~$m(m+1)$ switching points. 

In the structure of the Barabanov norm given by Theorem~\ref{th.40}, 
the set~$D$ is finite and hence, can be neglected: the integral over 
 this set is equal to zero. The sets of intervals~$P$ and~$Q$ are finite. 
 The formula~(\ref{eq.g}) becomes 
\begin{equation}\label{eq.g-f}
g(\gamma)\ = \ g(0) \ + \ 
\sum_{q_i \subset [0, \gamma]}\, \Delta(q_i) \ + \ \sum_{p_j \subset [0, \gamma]}\, \Delta(p_j)\, . 
\end{equation}
The existence of the  invariant sphere~$S$ with parameters~$\Delta(p_j) = s_j  \in [m_1, M_i], \, p_j \in P$, is equivalent to the condition 
\begin{equation}\label{eq.si-f}
\sum_{p_j \in P}\, s_j  \ = \ 
 - \ 
\sum_{q_i \in Q}\, \Delta(q_i) \, , 
\end{equation}
where all summations are finite. 
 
\bigskip

\begin{center}
\large{\textbf{7. Examples and applications }}
\end{center}
\bigskip

\begin{center}
\textbf{7.1. Stability of a matrix ball}
\end{center}
\bigskip

Consider the control 
set~$\cA(A_0, r) \, = \, \{A: \ \|A-A_0\|_F\, \le r\}$ being a  ball  
in the space of~$2\times 2$ matrices centered at a given matrix~$A_0$. 
Here~$\|X\| = \bigl({\rm tr}(X^TX)\bigr)^{1/2}\, = \, \bigl( \sum_{i,j} |X_{ij}|^2\bigr)^{1/2}$ is the Frobenius norm. 
Thus, one needs to decide the stability of the system~$\dot \bx = A(t)\bx$, where for 
each~$t$,~$A(t)$ is an arbitrary matrix at a distance of at most~$r$ from~$A_0$.  

 The matrix~$A_0$ is assumed to be 
Hurwitz, so, $\cA(A_0, 0)$ is stable.  Clearly, for every~$r>0$, the set $\cA$ is convex and irreducible. 
The remaining assumption that no matrix from~$\cA$ possesses 
a  nonnegative eigenvalue is solved by the following simple lemma. 
The norm of all vectors below is Euclidean. 
\begin{lemma}\label{l.100}
The ball~$\cA (A_0, r)$ contains a matrix with a nonnegative eigenvalue if and only if 
there exists~$\bx$ such that~$\|\bx\| = 1$ and the distance from~$A_0\bx$ to the 
ray~$\{\lambda \bx : \ \lambda \ge 0\}$ does not exceed~$r$. 
\end{lemma}
{\tt Proof}. Let~$A \in \cA(A_0, r)$ be such that~$A\bx = \lambda \bx, \, \lambda \ge 0$; then~$\|A_0\bx - \lambda \bx\| = \|(A-A_0)\bx\| \, \le \, \|A-A_0\|_F \|\bx\|\, \le \, r$.
Hence, the distance from~$A_0\bx$ to the ray~$\ell = \{\lambda \bx : \ \lambda \ge 0\}$ does not exceed~$r$. Conversely, if this  distance does not exceed~$r$, then 
taking~$\lambda \ge 0$ such that~$\|A_0\bx - \lambda \bx\| \le r$ 
and~$A = A_0 - \by\, \bx^T$, where~$\by = A_0\bx - \lambda \bx$, we see that
$\|A-A_0\|_F = \|\by\, \bx^T\|_F = \|\by\|\, \|\bx\| \le r$ 
and~$A\bx = (A_0 - \by\, \bx^T)\bx \, = \, A_0\bx - \by \, = \, \lambda \bx$.  

  {\hfill $\Box$}
\medskip

Note that the conditions of Lemma~\ref{l.100} are elementary verified
by denoting~$\bx = (\cos\, s\, , \, \sin \, s)^T$ and finding the 
minimal distance from~$A_0\bx$ to the ray~$\{\lambda \bx : \ \lambda \ge 0\}$
over all~$s$. The point of minimum~$s$ is evaluated in an explicit form. 
Let $\sigma_1 \le \sigma_2$ be the singular values of~$A_0$
(we hope that this standard notation will not be confused with the Lyapunov exponent).

\begin{lemma}\label{l.110}
The ball~$\cA (A_0, r)$ contains a degenerate matrix if and only if~$r$ is not less than 
the minimal singular value~$\sigma_1(A_0)$.
If~$\sigma_1(A_0) = r$, then a unique degenerate matrix 
in that ball is~$A_0 \, - \, r  \bu\,  \bv^T$ , where $\bu, \bv$ 
are unit eigenvectors of the matrices~$A_0^TA_0$ and~$A_0A_0^T$, respectively, 
corresponding to the eigenvalue~$r^2$. 
\end{lemma}
{\tt Proof}. If~$A\bx = 0, \, \|\bx\| = 1$, then~$\|A_0\bx\|  =  \|(A_0-A)\bx\| \le r$, hence, 
$\sigma_1(A_0) \le r$. If~$\sigma_1(A_0)  = r$, then 
the only unit vector~$\bx$ satisfying the inequality~$\|A_0\bx\| = r$
is the corresponding eigenvector of~$A_0A_0^T$. Then, as in the proof of Lemma~\ref{l.100}, 
we conclude that~$A\bx = 0$ and~$\|A-A_0\|\le r$ if and only if~$A= A_0 - r \bu\, \bv^T$.

  {\hfill $\Box$}
\medskip

\begin{figure}[h!]
\centering
 	{\includegraphics[scale = 0.6]{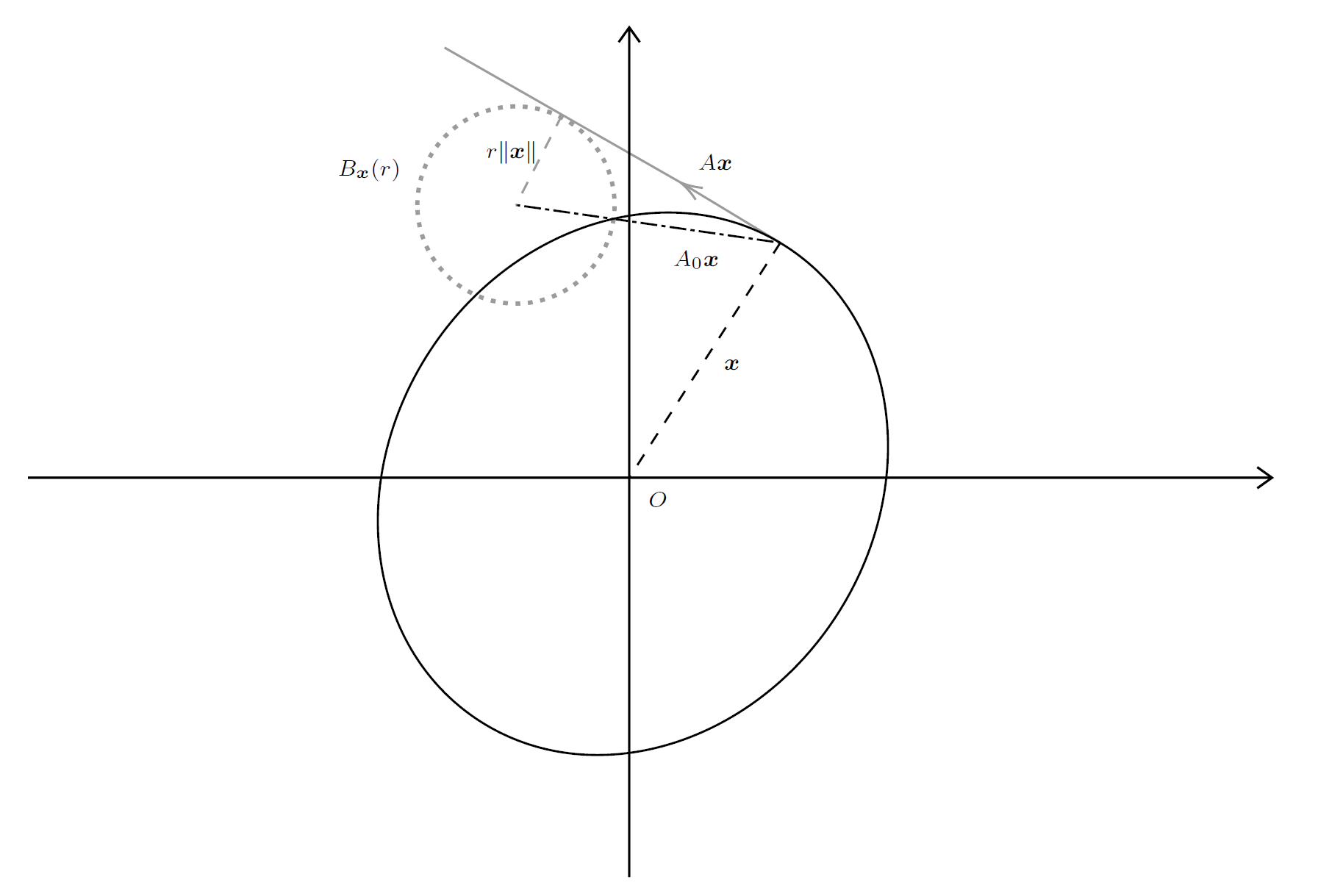}}
 	\caption{The Barabanov norm for a matrix ball. The vector~$A\bx$ touches~$B_{\bx}(r)$
 	and~$S$.}
 	\label{fig.r03475}
 \end{figure}

Now we can apply Theorem~\ref{th.30}. For each~$\bx$, the set~$\bx + A\bx, \,  A\in \cA$, 
is a disc of radius~$r\,\|\bx\|$ centered at the point~$\bx + A_0\bx$. Denote this disc by~$B_{\bx}(r)$, see Fig.~\ref{fig.r03475}. 
Then, both~$a_{\bx}$ and $b_{\bx}$ are tangents to~$B_{\bx}(r)$, provided that
they lie in the corresponding half-planes with respect to the line spanned by~$\bx$. 
If the angle between those tangents contains the origin, then~$a_{\bx}$ and $b_{\bx}$
are both well-defined, otherwise,  
only one of them is defined being the ``highest'' tangent, which makes the smallest angle with the vector~$\bx$. The corresponding left round leading trajectory is found by the 
ODE~$\dot \bx \, = \, A^{(l)}\bx$, 
with~$A^{(l)} \, = \, A_0 \, + \, \frac{1}{\|\bx\|}\, \bz \,\bx^T$, where~$\bz$ is the vector 
from the center of~$B_{\bx}(r)$ to the point of tangency.  In fact, one can 
use an arbitrary vector field~$F(\bx)$ which at every~$\bx$ is directed 
along~$a_{\bx}$; then the solution of ODE~$\dot \bx \, = \, A^{(l)}\bx$ describes the 
left leading trajectory. The same can be said on the right trajectory. 

To construct  the Barabanov norm, assume we again that~$\sigma(\cA) = 0$
after a proper shift,
where~$\cA = \cA(A_0, r)$, and invoke Lemma~\ref{l.110}. 
\begin{theorem}\label{th.80}
Consider the ball~$\cA = \cA(A_0, r)$ and assume that~$\sigma(\cA) = 0$; 
then~$r\le \sigma_1(A_0)$. If this inequality is strict, then all matrices from~$\cA$
are non-degenerate, the Barabanov norm is unique and described by the 
periodic leading left (or right) trajectory. 

If~$r = \sigma_1(A_0)$, 
then the only degenerate operator in~$\cA$
is~$A_0 - r  \bu\,  \bv^T$, where $\bu, \bv$ 
are unit eigenvectors of the matrices~$A_0^TA_0$ and~$A_0A_0^T$, respectively, 
corresponding to the eigenvalue~$r^2$, the Barabanov norm is 
described by the leading left (or right) trajectory from~$\bx(0) = \bv$ to
$\bx(+\infty) = - \bv$ and by its symmetric trajectory about the origin. 
\end{theorem}
{\tt Proof}. If~$r < \sigma_1(A_0)$, then~$\cA$ contains no degenerate matrices
(Lemma~\ref{l.110}, and~$S$ coincides with the leading periodic trajectory (Theorem~\ref{th.20}). If~$r =  \sigma_1(A_0)$, then~$\cA$ contains a unique degenerate matrix 
with  kernel spanned by the vector~$\bv$. Now the result follows from Lemma~\ref{l.70}
with $\bs = \bv$.

  {\hfill $\Box$}
\medskip

Applying now Theorem~\ref{th.50}, we obtain
\begin{cor}\label{c.50}
The Barabanov norm of a matrix ball is always unique. 
\end{cor}
\begin{remark}\label{r.130}
{\em It is not difficult to write an analytical expression for the
tangents to the ball~$B_{\bx}(r)$ from the point~$\bx$. This way we get an ODE, whose solution describes a unit sphere~$S$. From this equation one can conclude that~$S$ is in general not an ellipse, although can look similar to it (Fig.~\ref{fig.r03475}). }
\end{remark}
\bigskip

\begin{center}
\textbf{7.2. Stability of linear switching systems with a noisy data}
\end{center}
\bigskip 

Infinite control sets naturally appear when the matrices of the system are not given precisely but with possible deviation. At each moment~$t$ the deviation 
can be an arbitrary matrix from a given compact ``noise set''~$\cD$. 
Consider a finite  control set~$\cA = \{A_1, \ldots , A_m\}$; then the 
noisy linear switching system  is 
$$
\dot \bx \ = \ \bigl(A(t) + \Delta(t)\bigr)\bx\, , \qquad 
 A(t) \in \{A_1, \ldots , A_m\}, \, \Delta (t) \in \cD\, , \qquad  t \in [0, +\infty). 
 $$
If $m=1$ and~$\cD$ is a Frobenius ball, then 
we have the system from the previous subsection. 
The case of general~$m$ is considered similarly.  
In Theorem~\ref{th.80} the leading rays~$a_{\bx}, b_{\bx}$ are chosen from 
$2m$ tangents drown from the point~$\bx$ to the equal discs of radius~$r\, \|\bx\|$ 
centered at the points $\bx + A_i\bx, \ i=1, \ldots , m$.  

 If~$\cD$ is a matrix polyhedron, then 
so is~$\cA + \cD$, hence the problem is reduced to the 
control set with a finite number of matrices (vertices of this polytope). 
This is the case, for instance, if~$\cD$ is given by elementwise 
constraints:~$|\Delta_{ij}|  \le  \varepsilon_{ij}, \ i,j \in \{1,2\}$. 

For example,  $\cD$ can be a ball of radius~$r$ in the operator~$L_{\infty}$-norm in~$\re^2$, i.e., $\|\bx\|_{\infty} = \max\, \{|x_1|, |x_2|\}$, each set~$\bx + (A_i+ \cD)\, \bx$
is a square of radius~$r\,\|\bx\|$.

\begin{center}
\textbf{7.3. Numerical example}
\end{center}
\bigskip 

\begin{figure}[h!]
\centering
 	{\includegraphics[scale = 0.9]{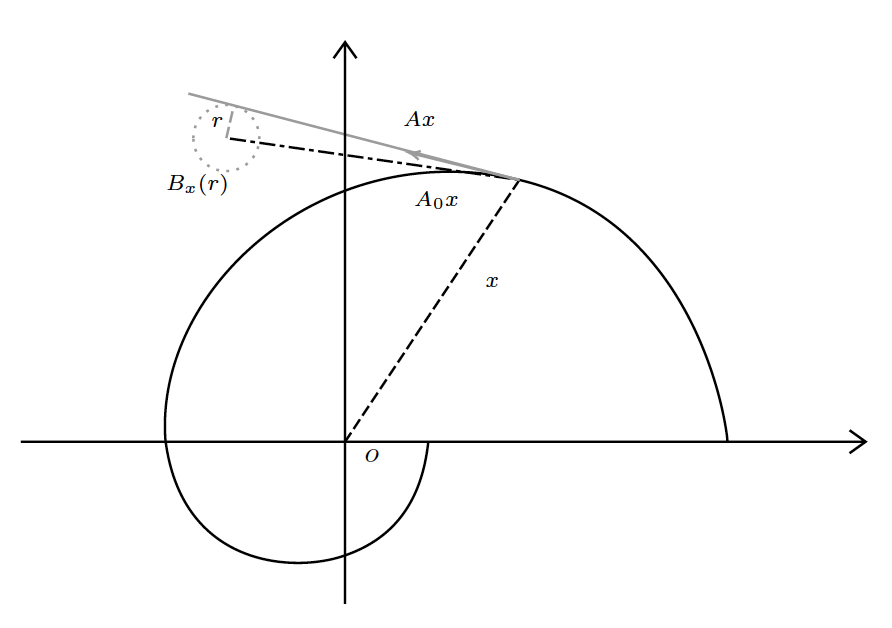}}
 	\caption{The leading trajectory for the matrix ball of radius~$r=0.1$}
 	\label{fig.r0}
 \end{figure}

The following matrix: 
$$
A_0\ = \ 
\left( 
\begin{array}{rr}
- 0.2 & - 1\\
1 &  - 0.5
\end{array}
\right)
$$ 
is Hurwitz, $\sigma(A_0) = -0.35$. We compute the radius~$R$ of the minimal 
unstable Frobenius ball with center~$A_0$.  
\begin{figure}[h!]
\centering
 	{\includegraphics[scale = 0.6]{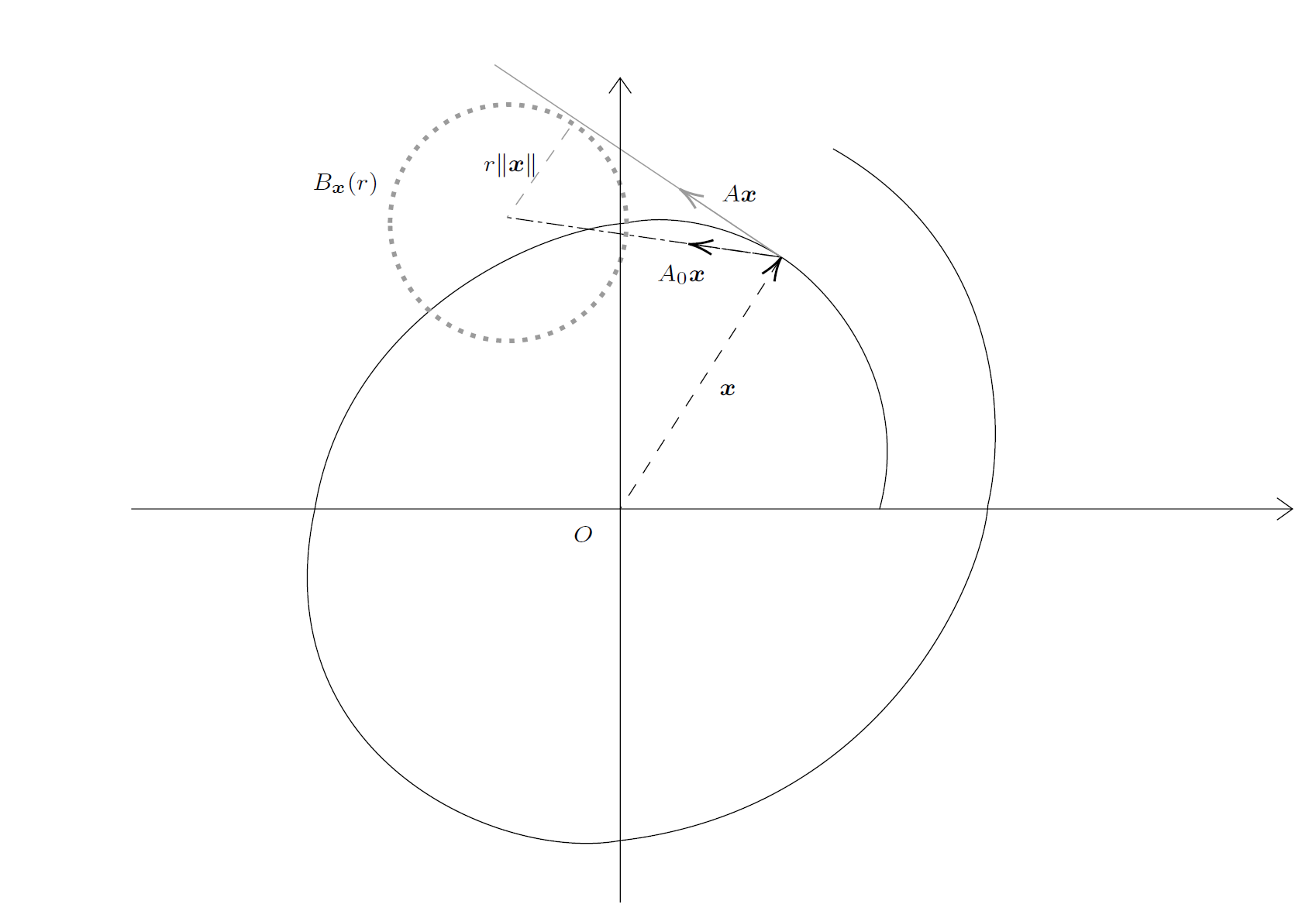}}
 	\caption{The leading trajectory for the matrix ball of radius~$r=0.4$}
 	\label{fig.r04}
 \end{figure} 
The computation is done 
by bisection in~$r$. For $r=0$, we have the one point ``ball'' $A_0$, 
which is clearly stable. For $r=0.4$ we construct left and right extremal trajectories. 
  Fig.~\ref{fig.r04} shows the left trajectory. By Theorem~\ref{th.80},  at every point~$\bx$, the leading vector~$A^{(l)}(t)\bx$ is directed along the tangent line 
from~$\bx$ to the ball~$B_{\bx}(r)$.
This trajectory is found from the differential 
equation~$\dot \bx(t) = A^{(l)}(t)\bx(t)$, we see that it 
 increases after one rotation. Hence, $\cA(A_0, 0.4)$ is unstable. Drawing the unique 
leading trajectory for~$r=0.1$ we see that it decreases, hence, 
by Theorem~\ref{th.30},  the ball  $\cA(A_0, 0.1)$ is stable (Fig.~\ref{fig.r0}). 
Thus,~$r\in [0.1\, , \, 0.4]$. Then  the bisection in~$r$  after 12 iterations gives the value~$0.3475...$ which approximates~$R$ with precision~$10^{-4}$.
 The left leading trajectory for the ball~$\cA(A_0, R)$ cycles and 
gives the sphere~$S$ of a unique Barabanov norm (Corollary~\ref{c.50}). 
It is drawn in Fig.~\ref{fig.r03475}.  
 \begin{remark}\label{r.120}
 {\em Note that~$S$ is a solution of a nonlinear ODE, and it is not ellipse 
 (Remark~\ref{r.130}). Also  observe that~$R = 0.3475... < 0.35 = |\sigma(A_0)|$. Thus, in this example,  
 the radius of the minimal  unstable ball is less than the modulus of the Lyapunov exponent of its center~$A_0$. }
 \end{remark}

\bigskip

\begin{center}
\large{\textbf{Appendix}}
\end{center}
\bigskip

{\tt Proof of Lemma~\ref{l.40}.} We have~$A^{(r)}\bx = - s A^{(l)}\bx$ with~$s>0$. 
In the basis~$\be_1 = \bx,\, \be_2 = A^{(l)}\bx$, we have 
$$
A^{(l)} \ = \ 
\left( 
\begin{array}{rr}
0& a\\
1 & b
\end{array}
\right) \ ; 
\qquad  A^{(r)}  \ = \ 
\left( 
\begin{array}{rr}
0& c\\
-s & d
\end{array}
\right)\ .  
$$
Consider the following matrix: 
$$
A_s\ = \ \frac{s}{s+1}\ A^{(l)} \ + \  \frac{1}{s+1}\ A^{(r)} \quad = \quad
\frac{1}{s+1}\, \left( 
\begin{array}{rr}
0& \ s\, a \, + \, c\\
0 & \  s\, b \, + \, d
\end{array}
\right)\, . 
$$
Since~$A_s$ is degenerate and belongs to~${\rm co}\, (\cA)$, it follows that 
$\sigma (\cA) \ge 0$, which proves the first claim. Now assume that $\sigma (\cA) = 0$. 
Then the eigenvalues of all matrices from~${\rm co}(\cA)$ have non-positive real parts.
If~$\, s\, a \, + \, c \, \ne \, 0$, then one can take a number~$s'$ close to~$s$ so that 
the determinant of the 
matrix~$\frac{s'}{s'+1}\ A^{(l)} \ + \  \frac{1}{s'+1}\ A^{(r)}$
is negative. Hence, one of its eigenvalues is positive, which is prohibited since this 
matrix belongs to~${\rm co}(\cA)$. 
Thus,~$s\, a \, + \, c \, = \, 0$, therefore, 
$$
A^{(l)} \ = \ 
\left( 
\begin{array}{rr}
0& a\\
1 & b
\end{array}
\right)\ ; 
\qquad  A^{(r)}  \ = \ 
\left( 
\begin{array}{rr}
0& -s\, a\\
-s & d
\end{array}
\right)\ .  
$$
Since~$A^{(l)}, A^{(r)}$ have eigenvalues with non-positive real parts, 
  their traces  are also non-positive, hence
$b\le 0, d\le 0$.  Consider now two cases for the coefficient~$a$. 
\smallskip 

1) $\ a\ne 0$. In this case~$a<0$, otherwise,~${\rm det}\, A^{(l)} < 0$, hence,
one of the eigenvalues of~$A^{(l)} < 0$ is positive. 
If $b=0$, then the eigenvalues of $A^{(l)}$
are~$\, \pm \, i \sqrt{|a|}$. Therefore, $A^{(l)}$ defines an ellipsoidal rotation and 
hence, $\cA$ has complex dominance (Proposition~\ref{p.10}). 
 The same happens if~$d=0$. 
Thus, it can be assumed that~$b<0, d< 0$. The pair $\{A^{(l)}, A^{(r)}\}$ 
is irreducible. Indeed, if they share a common real eigenvector, then 
the same eigenvector is associated to the matrix~$A_s$, whose eigenvectors are~$\be_1$ and
$\be_2$ (remember that~$s a + c =  0$). 
 None of them is an eigenvector of~$A^{(l)}$, because~$a\ne 0$. 

Thus, the pair $\{A^{(l)}, A^{(r)}\}$ is irreducible and its 
Lyapunov exponent is equal to zero, hence, it possesses a Barabanov norm. 
A unit sphere of the Barabanov 
norm of two $2\times 2$ matrices is composed with arcs of trajectories of the two 
equations~$\dot \bx \, = \,  A^{(l)}\bx$ and~$\dot \bx \, = \,  A^{(r)}\bx$;   
 and every switching point~$\bv$ between those arcs 
satisfy~$(A^{(r)} \, - \, \lambda A^{(l)})\bv\, = \, 0$, where~$\lambda$  
is found from the equation~${\rm det}\, (A^{(r)} \, - \, \lambda A^{(l)})\, = \, 
0$, see~\cite[Corollary~1]{Mu1}. 
 In our case this equation reads~$a(\lambda +s)^2\, = \, 0$. 
Hence,~$\lambda = - s$ and the solution~$\bv$ is unique up 
to normalization. Hence, if there are switching points, then 
they are~$\pm \bv$. We see that the arc of~$S$ connecting those points 
$\bv$ and~$-\bv$ satisfies one equation, say~$\dot \bx \, = \,  A^{(l)}\bx$. 
Therefore, $- \bv \, = \, e^{t\, A^{(l)}} \bv$ for some~$t>0$. 
This means that one of the eigenvalues of the matrix~$e^{t\, A^{(l)}}$ is~$-1$. 
Therefore,~$A^{(l)}$ has eigenvalues~$\pm \alpha i $ for some~$\alpha \ne 0$, hence, 
we arrive again at the case of complex dominance. 
\smallskip 

2) $\ a= 0$. In this case~$A^{(l)}, A^{(r)}$ are degenerate, have 
a common image parallel to~$\be_2$ and  kernels given by 
the equations~$x_1 + bx_2 = 0$ and $x_1 - \frac{d}{s}\, x_2 = 0$. 
These kernels are different since~$b < 0$ and~$\, - \frac{d}{s} \, > \, 0$. 
\smallskip 

Thus, $A^{(l)}$ and $A^{(r)}$ belong to one pencil. By Lemma~\ref{l.35},~$S$
contains a segment parallel to~$\be_2$ bounded by the kernels. 
Then~$\bx$ is its interior point, because it does not belong to the kernels.  

  {\hfill $\Box$}
\medskip

{\tt Proof of Proposition~\ref{p.50}}. If $A$ is not an extreme point in~${\rm co}(\cA)$, then 
$A = \frac12 \bigl(A_1 + A_2 \bigr)$ for some~$A_1, A_2 \in {\rm co}(\cA)$.
Let~$\bx \in {\rm Ker}(A)$. Then~$A_1\bx + A_2\bx = 0$. 
If~$A_1\bx \ne 0$, then the vectors~$A_1\bx$ and $A_2\bx$ 
have opposite directions. In this case, the angle between~$a_{\bx}$ and~$b_{\bx}$
is at least~$\pi$. It cannot be bigger than~$\pi$, the proof is the same as in Step~1 
of the proof of Theorem~\ref{th.20}. If it is equal to~$\pi$, then 
by Lemma~\ref{l.40}, $A$ belongs to some pencil~$\cP[A_1,A_2]$ and 
$\bx$ is an interior point of the reverse segment. This implies that 
the operators generating this pencil are extreme, i.e., belong to~$\cA$. 
If $A_1\bx = 0$, then we  consider the set of operators~$\cQ \subset  {\rm co}(\cA)$
with the same kernel spanned by~$\bx$. The  images of those operators are lines of 
support for~$S$ at the point~$\bx$, see~Proposition~\ref{p.60}. Let~$\ell_1, \ell_2$ be 
the extreme lines among them. Denote by~$\cQ_i$ the set of operators 
from~$\cQ$ with the image parallel to~$\ell_i,  \, i = 1, 2$. 
Each~$\cQ_i$ is a segment (in the space of matrices) of a line passing though the origin. 
Let~$A_i$ be the end of this segment most distant from the origin. 
Then, $A_i$ is an extreme point of~${\rm co}(\cA)$, i.e. $A_i \in \cA$, and 
$A$ is a linear combination of~$A_1, A_2$. 

  {\hfill $\Box$}
\medskip 

\textbf{Acknowledgements}. 
The research  is performed with the support of the Theoretical Physics and Mathematics Advancement Foundation ``BASIS'', grant no. 22-7-1-20-1.

 \end{document}